\documentclass[12pt]{article}

\usepackage[utf8]{inputenc}  
\usepackage{amsmath, amssymb}  
\usepackage{graphicx}  
\usepackage{hyperref}  
\usepackage{geometry}  
\usepackage{color}  
\usepackage{cite}  
\usepackage{authblk}  
\usepackage{booktabs}
\usepackage{tabularx}

\usepackage{lipsum}
\usepackage{amsfonts}
\usepackage{graphicx}
\usepackage{subfigure}
\usepackage{float}    
\usepackage{comment}
\usepackage{epstopdf}
\usepackage{algorithmic}
\usepackage{bbm}
\ifpdf
  \DeclareGraphicsExtensions{.eps,.pdf,.png,.jpg}
\else
  \DeclareGraphicsExtensions{.eps}
\fi

\newcommand{\R}{\mathbb{R}}

\newcommand{\x}{\mathbf{x}}
\newcommand{\X}{\mathbf{X}}

\newcommand{\uv}{\mathbf{u}}
\newcommand{\vuv}{\mathbf{v}}

\newcommand{\UV}{\mathbf{V}}
\newcommand{\y}{\textbf{y}}
\newcommand{\eps}{\varepsilon}

\newtheorem{defi}{Definition}[section] 

\title{Numerical study on a multi-dimensional pressureless Euler-type model with non-local interactions and chemotaxis for collective cell migration}

\author[1]{M. Menci}
\author[2]{R. Natalini}
\author[3]{T. Tenna}

\affil[1]{Università Campus Bio-Medico di Roma, Rome, Italy}
\affil[2]{Istituto per le Applicazioni del Calcolo, Consiglio Nazionale delle Ricerche, Rome, Italy}
\affil[3]{Université Côte d'Azur, Nice, France - Sapienza Università di Roma, Rome, Italy}

\hypersetup{
    pdftitle={Numerical study on a multi-dimensional pressureless Euler-type model with non-local interactions and chemotaxis for collective cell migration},
    pdfauthor={M. Menci, R. Natalini, T. Tenna}
}

\begin{document}

\date{}
\maketitle

\begin{abstract}
In this paper we propose a numerical study of macroscopic models for collective cell migration, focusing on a multi-dimensional pressureless Euler-type model with non-local interactions coupled with chemotaxis, rigorously derived from microscopic dynamics. Different mechanical interactions are investigated, including attraction-repulsion effects. Moreover, the model is extended to the case of different populations of interacting cells. The validity of such macroscopic model and its agreement with the microscopic dynamics is finally assessed through a parameter estimation analysis in a specific setting.
\end{abstract}

\section{Introduction}

Many different real-world collective phenomena exhibit self-organization behaviors, usually driven by long-range mechanical interactions among agents. In the biological field, collective cell migration is also induced by chemical signals, detected by the cells and guiding their oriented motion. Several processes are regulated by the response of cell populations to the change of chemical concentrations in the environment and different mathematical models
were developed for an accurate description of these phenomena \cite{murray2007mathematical, perthame2006transport}. Starting from the seminal Keller-Segel model introduced in \cite{keller1970initiation,keller1971model}, 
where the evolution of the density of cells is described by a parabolic equation and the concentration of a chemoattractant is given by a parabolic (or elliptic) equation, a wide number of variations have been proposed. 
Among others, in \cite{nanjundiah1973chemotaxis} authors present a minimal version of the original model, that has been analytically investigated in different directions. 
A more complete analysis and some variations can be found in \cite{hillen2009user, painter2019mathematical}.
In more recent years, hyperbolic models have been developed. In particular, the models in \cite{PhysRevLett.90.118101, tosin2006mechanics, serini2003modeling} 
have been proposed to describe the early stages of the vasculogenesis process and they consist of a quasi-linear hyperbolic system, directly derived from macroscopic balance laws, coupled with a parabolic equation for the chemoattractant. The advantage of the hyperbolic approach over parabolic Keller–Segel-type models relies in its ability to reproduce emerging pattern, i.e. the formation of the network of capillaries, as numerically observed in \cite{natalini2015numerical}.
Moreover, the models account for the presence of a pressure term as a key factor, which authors consider to prevent overcrowding of cells, on a phenomenological basis.

Nevertheless, the behavior of cells at the microscopic scale seems to play a crucial role in modeling collective cell migration, providing a detailed description of the mechanical interactions among cells. In particular, second-order differential models allow to take into account the crucial role of persistence effects. This key feature of cells behavior is not observed in other collective phenomena, for instance in pedestrian dynamics, where accelerations are almost instantaneous and inertia-based effects are negligible.
In recent years, a novel class of \textit{hybrid coupled models} has been proposed, gathering the advantages of a multiscale description \cite{di2015hybrid, costanzo2020hybrid, di2018discrete}. These models are based on a \textit{discrete in continuous} approach: cells are described as discrete entities moving in the environment and chemotactic signals, influencing the dynamics, as continuous concentrations.  Microscopic and macroscopic scales are strongly coupled: the dynamics of cells are influenced by the gradient of the chemical concentration and, at the same time, the evolution of the chemical signal is influenced by cells positions.
From the analytical point of view, well-posedeness results for generalized versions of hybrid models have been investigated over the last years \cite{menci2022existence, menci2019global, menci2023existence, menci2023coupled}.

From the modeling and numerical perspective, the detailed level of description of the microscopic approach can become unsuitable if the phenomena under investigation involve high number of cells. A macroscopic approach, modeling the evolution in time of macroscopic quantities of interest, as cell density and velocity field, better fulfills the requirements. For coupled hybrid systems, the issue concerning the macroscopic counterpart of the microscopic dynamics naturally arises.
In this perspective, the recent results in \cite{natalini2021mean} represent the first step toward rigorous derivation of a macroscopic model.
Here authors rigorously proved the mean-field limit of the microscopic model towards a kinetic equation coupled with chemotaxis, and its hydrodynamic limit, under strong assumptions. 
The main feature of the resulting Euler-type model is the absence of pressure terms, in the face of non-local integral terms which keep into account the microscopic interactions, establishing a clear link between the macroscopic model and the microscopic dynamics. As a drawback, the absence of the pressure term in Euler-type systems leads to possible blow-up of the solutions \cite{breniergrenier1998, bouchutjames1999, carrillo2016critical}.
Non-local pressureless Euler systems have been derived and investigated in different directions \cite{natalini2020mean, motsch2011new, figalli2018rigorous, albi2018pressureless}.
In particular, a wide number of results deals with alignment kind of interactions, inspired by the seminal microscopic alignment model in \cite{cucker2007emergent}.

\paragraph{Aim of the paper}
In this paper we focus on the macroscopic scale of description, in a multi-dimensional setting. A numerical investigation of different models is performed, comparing seminal models of the literature with the recent pressureless Euler-type model coupled with chemotaxis derived in \cite{natalini2021mean}.
The absence of pressure enables the development of a model that does not require phenomenological terms, relying instead on nonlocal terms that retain information about microscopic interactions. In this perspective, we extend the kind of mechanical interactions included in the model, considering attraction-repulsion effects 
together with alignment. The non-locality of interactions and the coupling with the chemical signal represent the main features of the model with respect to previous works of the literature.
Moreover, the original structure of the dynamics is extended to a multi-population case, accounting for different types of interacting cells. This choice allows the construction of more complex and realistic models, where nonlocal interactions are characteristic of each individual species.
Inspired by successful applications of microscopic hybrid models in simulating Cancer-On-Chip dynamics \cite{bretti2021estimation, bretti2024}, 
we test the applicability of the macroscopic model in this context through a calibration procedure. Due to the current lack of experimental data on the macroscopic scale, we present an estimation procedure linking the microscopic dynamics, already calibrated on real experimental data, and the macroscopic one. This strategy enables to generate synthetic data from the microscopic model, in order to obtain a greater number of different scenarios to assess the macroscopic one.

\paragraph{Structure of the paper}
The rest of the paper is organized as follows. In Section \ref{sec:models} we analyze in details the coupled models, starting from the microscopic scale to the macroscopic one. In particular, we present the extension of the non-local Euler model to the case of two populations of cells and including different interactions.
The numerical study of the models is deferred to Section \ref{Section: Numerical_Simulations}, aiming at investigating and comparing the role of the different terms. Section \ref{sec:estimation} is devoted to parameter estimation for the model introduced in Section \ref{sec:twopop_model}. All the details about the numerical schemes adopted for the approximation of the solutions to the different models are reported in Appendix \ref{Appendix_Numerical_Scheme}.


\section{Description of the mathematical models}\label{sec:models}

In this section we delve in details the mathematical models for chemotactic cell migration considered in our study, highlighting the role of the different terms. 
In the following, we distinguish between microscopic and macroscopic coupled models, referring to the scale of observation considered for cell dynamics.
In all the considered models, the evolution of chemoattractant is modeled by a linear diffusion equation with a source term which depends on the positions of cells (for the microscopic one) or the density of cells (for the macroscopic ones). We start by reviewing existing models from the literature, then incorporate additional nonlocal effects, both for the microscopic and the macroscopic scale. Finally, we extend the model to account for multiple interacting populations, focusing on the dynamics between immune and cancer cells at the macroscopic level. 

\subsection{Microscopic models}
In all the works of the literature \cite{di2015hybrid, di2018discrete, costanzo2020hybrid}, cells are modeled as circles, endowed with a radius $R>0$, whose dynamics is described by time update of positions and velocities of their centers.
From a mathematical point of view, the general structure of a hybrid model consists of a system of second order ordinary differential equations, modeling the evolution in time of position $\textbf{x}_i \in \mathbb{R}^d$ and velocity $\textbf{v}_i \in \mathbb{R}^d$ of the $i$-th cell, $i=1,...,N$, coupled with a reaction-diffusion parabolic partial differential equation which describes the evolution of the concentration of the chemoattractant $\varphi=\varphi(x,t)$. The model reads as:
\begin{equation}
\begin{cases}
    \label{eq:Hybrid_general}
    \displaystyle
    \displaystyle \dot{\x}_i=\vuv_i,\\ 
    \dot{\vuv}_i= \textbf{F}_i\left(t,\textbf{X}(t), \textbf{V}(t), \nabla \varphi(t,\textbf{x}_i) \right),\\[5pt]
	\partial_s \varphi=D \Delta \varphi -\kappa \varphi + g(x,\textbf{X}(s)), \ \quad \ s \in [0,t],
\end{cases}
\end{equation}
where $\textbf{X}=(\textbf{x}_1,\dots,\textbf{x}_N)\in \mathbb{R}^{Nd}$, $\textbf{V}=(\textbf{v}_1,\dots,\textbf{v}_N)\in \mathbb{R}^{Nd}$ denote the collections of position and velocity of each cell at any time in the time interval considered, $t \in [0,T]$.
The chemotaxis equation accounts for constant diffusion and degradation coefficients, respectively denoted with $D,\,\kappa >0$. The source term $g$ models the fact that the chemical signal is produced by the cells themselves through the dependence on the cells positions.  
In \cite{ di2018discrete,di2015hybrid,menci2022existence, menci2023microscopic} $g$ is defined as the sum of characteristic functions over circular domains of fixed radius $R>0$ modeling cells.
Throughout the rest of the paper, we keep the space variables in math mode style for simplicity of notation.
Functions $\textbf{F}_i$ model mechanical and chemotactic interactions influencing the dynamics of each cell $i=1,...,N$. In particular, in \cite{costanzo2020hybrid, natalini2021mean} $\mathbf{F}_i$ are of the following form:
\begin{equation}
\label{eq:F_i_expr}
    \mathbf{F}_i \left(t, \X(t), \UV(t), \nabla \varphi(t, \x_i) \right) = \frac{1}{N} \sum_{j=1}^N \gamma(\vuv_i-\vuv_j,\x_i-\x_j) + \eta \nabla \varphi(t,\x_i)+F_{ext}(\x_i, \vuv_i),
\end{equation}
where $\gamma: \mathbb{R}^d \times \mathbb{R}^d \to \mathbb{R}^d $  is the interaction function, $\varphi:[0,T] \times \mathbb{R}^d \to \mathbb{R}^d$ is the concentration of the chemoattractant and $F_{ext}$ models external forces.\\
Initial data are given by initial position and velocity of each particle, and initial concentration of chemoattractant
\begin{equation}
\label{InitialCondition_Hybrid}
\textbf{X}(0)=:\textbf{X}_0, \quad \textbf{V}(0)=:\textbf{V}_0, \quad \varphi(x,0)=:\varphi_0.
\end{equation}
Starting from the hybrid coupled system in \eqref{eq:Hybrid_general}, different kinds of cell-cell interactions have been considered, depending on the specific applications.
In the literature, the great majority of agent-based models accounts for alignment and/or attraction-repulsion kind of interactions.
In order to describe different contributions in \eqref{eq:F_i_expr}, we rewrite the interaction function $\gamma$ as the sum of two contributions, $\gamma=\gamma_{1}+\gamma_{2}$.
Among alignment models, the seminal Cucker-Smale model has been conceived to  model consensus in flocks of birds and animal groups \cite{cucker2007emergent, cucker2007mathematics, couzin2005effective}.
The force acting on each individual is a weighted average
of the differences of its velocity with the ones of each other member of the group, and it is proportional to the inverse of the distance between agents.
The interaction function is hence given in this form:
\begin{equation}
\label{eq:gamma_1}
\gamma_{1}(\vuv_i-\vuv_j, \x_i-\x_j) = \gamma_D(|| \x_i-\x_j ||) \left(\vuv_j-\vuv_i\right),
\end{equation}
where $\gamma_D$ is a decreasing function of the distance $|| \x_i -\x_j||$ between two agents.
In the following we assume $\gamma_D$ as the classical Cucker-Smale communication rate
\begin{equation}
    \gamma_D (||\x_i-\x_j||) = \displaystyle \frac{\beta}{\left(1+||\x_j-\x_i||^2\right)^\varsigma},
\end{equation}
where $\beta \in \R^+$ and $\varsigma \in \R^+$ is the Cucker-Smale exponent. Concerning attraction and repulsion effects, interactions are often modeled by interaction potentials \cite{d2006self, carrillo2013new} or introducing specific functions of the distance, modeling short-range repulsion and long-range attraction. 
In particular, in \cite{di2015hybrid, di2018discrete,szabo2006phase} authors assume an attractive linear elastic force, and a repulsion force proportional to the inverse of the distance between any two agents, defining a piecewise linear function:
\begin{equation}
\label{eq:gamma_2}
	\gamma_{2} (\x_i-\x_j) = 
	\begin{cases}
		\begin{aligned}
			-\omega_{rep} \left(\frac{1}{||\x_i-\x_j||}-\frac{1}{R_{rep}} \right) \frac{\x_j-\x_i}{||\x_i-\x_j||}, \quad \text{if } ||\x_i-\x_j|| &\leq R_{rep}, \\[7pt]
		  \omega_{adh} \, \left(||\x_i-\x_j||-R_{rep} \right) \frac{\x_j-\x_i}{||\x_i-\x_j||}, \quad \text{if } R_{rep}<||\x_i-\x_j|| &\leq R_{adh},
		\end{aligned}
	\end{cases}		
\end{equation}
with $\gamma_2 \equiv 0$ if $||\x_i-\x_j|| > R_{adh}$.
Here $\omega_{rep},\,\omega_{adh},\, R_{rep},\,R_{adh}$ are constant positive parameters. 
In particular repulsion occurs when the distance between two agents is less than $R_{rep}>0$, whereas adhesion occurs at a distance greater than $R_{rep}$ and less than $R_{adh}>R_{rep}$.
Several extensions have been proposed, generalizing \eqref{eq:Hybrid_general}. In \cite{di2015hybrid}, a hybrid system with non-local concentration is introduced, accounting for the fact that the gradient of chemical concentration influences cells not only in their centers, but also in their surroundings. The non-locality of the influence of chemical concentration leads to a more complex structure, which have been rigorously analyzed from a theoretical view point in \cite{menci2023coupled}.

\subsection{Macroscopic models}
\label{Macro_Models}
The quasi-linear hyperbolic model for cell migration proposed in \cite{PhysRevLett.90.118101} was originally conceived to describe vasculogenesis processes. The structure of the model is the following:
\begin{equation}
	\begin{cases}
 \label{eq:GambaPreziosi_intro}
		&\partial_t \rho + \nabla \cdot (\rho \uv)=0,\\
		&\partial_t (\rho\uv)+ \nabla \cdot (\rho \uv \otimes \uv) = \eta \rho \nabla \varphi  - \nabla P(\rho) - \alpha \rho \uv,\\
		&\partial_t \varphi=D\Delta \varphi  - \kappa \varphi + g(\rho),
	\end{cases}
\end{equation}
where $\eta,\,\alpha,\,D,\,\kappa$ are positive constants. The term $\alpha \rho \uv$ represents a damping term, which models dissipative interactions with the external environment. In the non-local macroscopic model \eqref{eq:NP_Eulero}, this force could be included in $F_{ext}$.
The function $g$ models the production rate of the chemoattractant.
In \cite{PhysRevLett.90.118101} a constant production is assumed, i.e. $g(\rho)= a\,\rho$ with $\rho>0$. In \cite{tosin2006mechanics} authors choose
$g(\rho)=\alpha_1 \rho/(1+\alpha_2 \rho^2)$, with $\alpha_1,\, \alpha_2 >0$, modeling the fact that, in aggregate state where cells are packed and density is high, the chemoattractant production almost vanishes.
The first equation is the classical mass conservation, since cells do not undergo mitosis or apoptosis during the considered phenomenon, hence the total number of cells is constant in time.  
The right-hand side of the second equation, modeling the different effects that influence changes in cell directional motion, is the sum of the following different contributions:
\begin{enumerate}
\item a chemotactic body force 
\begin{equation}
	f_{\text{chem}}=\eta \rho \nabla \varphi,
\end{equation}
where $\eta>0$ measures the strength of the cells response to concentration gradients of the chemical substance.
The linear dependence on $\rho$ corresponds to the assumption that each cell experiences a similar chemotactic action. In \cite{scianna2013} a saturation effect on the amount of chemoattractant is also included.  

\item a dissipative interaction with the substrate 
\begin{equation}
	f_{\text{diss}}=-\alpha \rho \uv.
\end{equation}
A linear dependence
of $f_{\text{diss}}$ on $\rho$ is assumed, and its strength is tuned through the constant parameter $\alpha>0$. 

\item the main feature relies on the pressure term 
\begin{equation}
	f_\text{press} = - \nabla P(\rho),
\end{equation}
which has been introduced with the aim of preventing overcrowding of cells.
In this sense $f_\text{press}$ is a phenomenological term modeling no-collision dynamics of cells.
The main issue is indeed the choice of the function $P$, which is not justified by empirical evidence.
In \cite{PhysRevLett.90.118101,
cavalli2007approximation,
tosin2006mechanics} the structure is inspired to the pressure law for isentropic gases, i.e. 
$P(\rho) \propto \rho^\Gamma, \Gamma >1$.\\

\end{enumerate}

Let us now focus on the macroscopic model derived in \cite{natalini2021mean}, which describes the evolution of cell density $\rho=\rho(x,t)$ and velocity field $\textbf{u}=\textbf{u}(x,t)$ through the following non-local pressureless Euler-type system coupled with chemotaxis
\begin{equation}
\begin{cases}
\label{eq:NP_Eulero}
	\partial_t \rho+ \nabla \cdot (\rho \uv)=0,\\
	\partial_t (\rho \uv) + \nabla \cdot (\rho \uv \otimes \uv) = \rho 
   \mathcal{I}
 + \eta \rho \nabla \varphi + \rho F_{ext},\\
	\partial_t \varphi = D \Delta \varphi - \kappa \varphi + g(\rho).
\end{cases}	
\end{equation}
Here $\eta,\,D,\,\kappa$ are positive constants and $\mathcal{I}=\mathcal{I}(x,t)$ denotes a non-local integral term, which in \cite{natalini2021mean} is defined as
\begin{equation}\label{eq:int_NP_Eulero}
\mathcal{I}(x)=\displaystyle 
\int_{\mathbb{R}^d} \gamma (\cdot-y, \uv(\cdot)-\uv(y))\rho(y)\, dy.
\end{equation}

The main difference with respect to \eqref{eq:GambaPreziosi_intro} is the absence of the pressure term, formally replaced by the non-local integral term $\mathcal{I}$ in \eqref{eq:int_NP_Eulero}. From a mathematical point of view, it is well-known that the presence of pressure may lead to some a priori estimates, preventing blow up of solutions and resulting as a stabilizing force \cite{kowalczyk2004stability}. 
The key factor of \eqref{eq:NP_Eulero} relies on the non-local integral term \eqref{eq:int_NP_Eulero}, which highlights the mechanical kind of interactions acting in the microscopic dynamics through the function $\gamma$.
Concerning an Euler-alignment system with the interaction function inspired by the Cucker-Smale model, the integral term can be rewritten as

\begin{equation}\label{int_align}
\mathcal{I}_1(x)=\displaystyle \int_{\mathbb{R}^d} 
 \frac{\beta}{\left(1+||x-y||^2\right)^{\varsigma}} \left( \uv(y)-\uv(x)\right) \rho(y)\, dy.
\end{equation}

Changing the integral function, one can simulate different kind of interactions. In particular, we will consider repulsion-attraction among cells, using the interaction function $\gamma_2$ in \eqref{eq:gamma_2}.
In this case, the integral term can be written as

\begin{equation}\label{int_attr_rep}
\mathcal{I}_2(x)=\displaystyle \int_{\mathbb{R}^d} 
 \gamma_2(x-y)\rho(y)\, dy.
\end{equation}

In Section \ref{Section: Numerical_Simulations} we present numerical simulations of the two macroscopic models here considered, in order to compare their specific features. The rigorous derivation of \eqref{eq:NP_Eulero} from the microscopic dynamics \eqref{eq:Hybrid_general} can be found in \cite{natalini2021mean}.

\subsection{Extension to two-populations models}
\label{sec:twopop_model}
In a microscopic description of the dynamics, the hybrid model \eqref{eq:Hybrid_general} has been extended to the case of two populations of interacting cells. In \cite{bretti2021estimation}, authors presented a hybrid model describing immune cells migration under the effects of chemicals released by cancer cells, which are fixed in the domain. The cell-cell interactions for the immune system dynamics is described by \eqref{eq:gamma_1}-\eqref{eq:gamma_2}. The dynamics of immune cells towards cancer cells is described by a radial repulsion term, depending on cells position. Let $\y_j$ be the position of the $j$-th cancer cell, then function $\gamma_3$ is introduced as
\begin{equation}
    \label{eq:gamma_3}
    \gamma_{3} (\x_i-\y_j) := 
		-\omega_{rep}^{tum} \left(\frac{1}{||\x_i-\y_j||}-\frac{1}{R_{rep}^{tum}} \right) \frac{\y_j-\x_i}{||\x_i-\y_j||}, \quad \text{if } ||\x_i-\y_j||\leq R_{rep}^{tum}.
\end{equation}
and $\gamma_3 \equiv 0$ if $||\x_i-\y_j|| > R_{rep}^{tum}$, where $\omega_{rep}^{tum}$. Here $R_{rep}^{tum}$ are constant positive parameters. With this choice, repulsion occurs when the distance between the $i$-th immune cell and the $j$-th cancer cell is less than $R_{rep}^{tum}$.
The final model for the description of immune cells migration in presence of cancer cells can be written as \cite{bretti2021estimation}:
\begin{equation}
    \begin{cases}
    \begin{aligned}
	\label{eq:FinalHybridModel_Tumor}
    \displaystyle
	 \displaystyle \dot{\x}_i=\vuv_i,\\
     \dot{\vuv}_i= \frac{1}{N} &\sum_{j=1}^N \left(\gamma_1(\vuv_i-\vuv_j,\x_i-\x_j) + \gamma_2(\x_i-\x_j) \right)\\ + &\sum_{j=1}^M \gamma_3(\x_i-\y_j) + \eta \nabla \varphi(t,\x_i)+F_{ext}(\x_i),\\
	\displaystyle \partial_s \varphi = D& \Delta \varphi -\kappa \varphi + \sum_{j=1}^M \chi (x-\y_j), \qquad s \in [0,t],
    \end{aligned}
\end{cases}
\end{equation}
 where $x \in \R^d$, $N$ is the total number of immune cells, $M$ is the total number of cancer cells of radius $R_{tum}$ and $\chi \in \mathcal{C}^1_c$ is a generic function, usually chosen as a weighted regularization of the characteristic function
 \begin{equation*}
     \mathbbm{1}(z) = \begin{cases}
         1 \quad &\text{if } \, ||z||< R_{tum},\\
         0 \quad &\text{otherwise}.
     \end{cases}
 \end{equation*}
System \eqref{eq:FinalHybridModel_Tumor} is complemented by initial conditions, given by microscopic positions and velocities of each immune cell:
\begin{equation}
    \x_i(0)=\x_i^0, \qquad \vuv_i(0)=\vuv_i^0.
\end{equation}

Since cancer cells do not move in the considered setting, the initial condition for the parabolic equation is given by the chemoattractant $\varphi$ produced at time $t=0$, namely
\begin{equation}
\label{initial_condition_tumor_chemoattractant}
    \varphi(x, 0) = \sum_{j=1}^M \chi (x-\y_j).
\end{equation}

Inspired by the structure of \eqref{eq:NP_Eulero} and its link with the microscopic scale previously described, we here propose a macroscopic Euler-type model keeping memory of the key factor of \eqref{eq:FinalHybridModel_Tumor}.
More specifically, the source term of the parabolic equation for the chemoattractant here depends on the density of cancer cells, as in \eqref{eq:FinalHybridModel_Tumor}. 
We denote by $\rho(x,t)$ and $\zeta(x,t)$ respectively the density of immune and cancer cells, respectively, at time $t$ in $x$. To model interactions between immune and cancer cells in \eqref{eq:NP_Eulero_Tumor}, we introduce a non-local integral term depending on the microscopic interaction modeled by $\gamma_3$ in \eqref{eq:NP_Eulero}, namely
\begin{equation}\label{int_rep_tumorali}
\mathcal{I}_3(x)=\displaystyle \int_{\mathbb{R}^d} 
 \gamma_3(x-y)\zeta(y)\, dy.
\end{equation}
This term links the microscopic interactions among cells with the macroscopic density evolution. The macroscopic two populations Euler-type model reads as
\begin{equation}
\label{eq:NP_Eulero_Tumor}
    \begin{cases}
        \partial_t \rho+ \nabla \cdot (\rho \uv)=0,\\
        \partial_t (\rho \uv) + \nabla \cdot (\rho \uv \otimes \uv) = \rho \mathcal{I} + \eta \rho \nabla \varphi + \rho F_{ext},\\
        \partial_t \varphi = D \Delta \varphi - \kappa \varphi + \chi * \zeta,
    \end{cases}
\end{equation}
where the non-local integral term $\mathcal{I}$ is  the sum of the alignment and attraction-repulsion contributions, 
\begin{equation}
    \mathcal{I} = \mathcal{I}_1 + \mathcal{I}_2 + \mathcal{I}_3,
\end{equation}
with $\mathcal{I}_1$, $\mathcal{I}_2$ and $\mathcal{I}_3$ respectively defined in \eqref{int_align}, \eqref{int_attr_rep} and \eqref{int_rep_tumorali} and cancer cells density $\zeta$ is assumed to be constant in time.

\section{Macroscopic coupled models: 2D numerical simulations} 
\label{Section: Numerical_Simulations}

In this section we present 2D numerical solutions to the models introduced in the previous sections.
The choice of the two-dimensional domain lies in the fact that most of the experimental data at our disposal refers to bi-dimensional settings. Indeed, biological experiments on two-dimensional devices offer a cost-effective tool to study cells interactions in a microfluidic environment.
To the best of our knowledge, up to now only 1D numerical simulations of the model in \eqref{eq:NP_Eulero} can be found in the literature \cite{menci2023microscopic}, and accounting only for alignment kind of interactions.

We approximate the solutions through finite difference schemes on a square spatial domain $\Omega=[0,1] \times [0,1]$, with a mesh size $\Delta x$ in $x$-direction and $\Delta y$ in $y$-direction. In all simulations we fix $\Delta x = \Delta y = 0.02$. The time interval is set to $[0,T]$, with $T=1$ and the numerical scheme adopt a time-variable discretization step $\Delta t$, according to the CFL condition. The numerical scheme is based on discrete kinetic approximations of the models, in which the transport term and the relaxation term are treated using a splitting technique \cite{aregba2000discrete, natalini1998discrete}. More precisely, the finite difference scheme for the transport term involves flux limiters in space, improving the numerical scheme conceived in \cite{menci2023microscopic} for the 1D case. The numerical approximation of the parabolic equation for the chemotaxis is based on a standard centered scheme. 
A detailed description of the numerical scheme is given in Appendix \ref{Appendix_Numerical_Scheme}.

The initial macroscopic density is obtained as the sum of $N$ Gaussian bumps whose amplitude is assumed to be of the order of a positive parameter $\sigma$. Heuristically, in the context of cell migration phenomena, $\sigma$ represents the non-dimensional average cell radius: each Gaussian bump, centered in random positions $\x_j$, $j=1,\dots, N$, approximates a cell endowed with its radius $\sigma$. Thus, the initial density is given by 
\begin{equation}
	\rho(x,0)=\frac{\tau}{2 \pi \sigma^2} \displaystyle \sum_{j=1}^N \exp\Biggl(-\frac{|x-\x_j|}{2 \sigma^2} \Biggr),
\end{equation}
where $\tau$ is a positive normalization parameter. We set $N=200$ and $\sigma=0.015$ to obtain the initial density shown in Fig. \ref{fig:density_initialcondition}. 
In the following tests we assume the same initial condition, setting null initial conditions for the velocity and the chemoattractant
\begin{equation}
    \rho_0(x)=\rho(x,0) \qquad \uv(x,0)=\uv_0(x) = 0, \qquad \varphi(x,0) = \varphi_0(x) = 0 \qquad \forall x \in \Omega.
\end{equation}

We assume homogeneous Neumann boundary conditions for the density $\rho$ and for the chemoattractant $\varphi$, and zero boundary conditions for the normal component of the velocity field $\uv$:
\begin{equation}
    \nabla \rho \cdot \mathbf{n} \, \big |_{\partial \Omega} = 0, \qquad \nabla \varphi \cdot \mathbf{n} \, \big|_{\partial \Omega} = 0, \qquad
        \uv \cdot \mathbf{n} \, \big|_{\partial \Omega} = 0. 
\end{equation}

We present different representative numerical simulations of the models in  \eqref{eq:NP_Eulero} and \eqref{eq:GambaPreziosi_intro}, showing the role of different effects involved in the dynamics. 

In all of the following numerical simulations we fix $\alpha=0$, $\kappa=1$, $D=0.1$, $\alpha_1=30$, $\alpha_2=0.2$, unless otherwise mentioned.

\subsection{Pressureless Euler with chemotaxis without non-local interactions}
We start our analysis simulating \eqref{eq:NP_Eulero} neglecting the non-local interaction term $\mathcal{I}$ together with external forces. The structure of the system hence reduces to a pressureless Euler system, coupled with chemotaxis. In our simulations, the finite-time blow-up which characterizes pressureless dynamics is observed, coherently with previous works of the literature on pressureless Euler (without chemotactic interactions) \cite{bouchutjinli2003, albi2018pressureless, carrillo2017review}.

\begin{figure}[H]
\centering
\subfigure[$\eta=0.1$]{\includegraphics[scale=0.3]{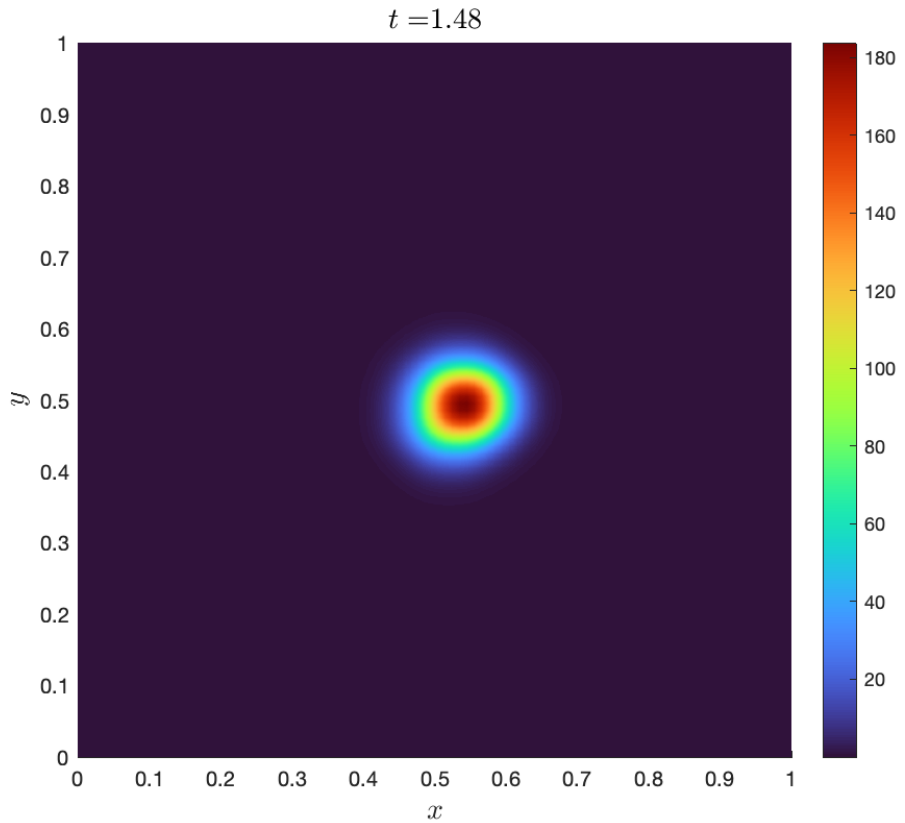}}
\subfigure[$\eta=0.5$]{\includegraphics[scale=0.3]{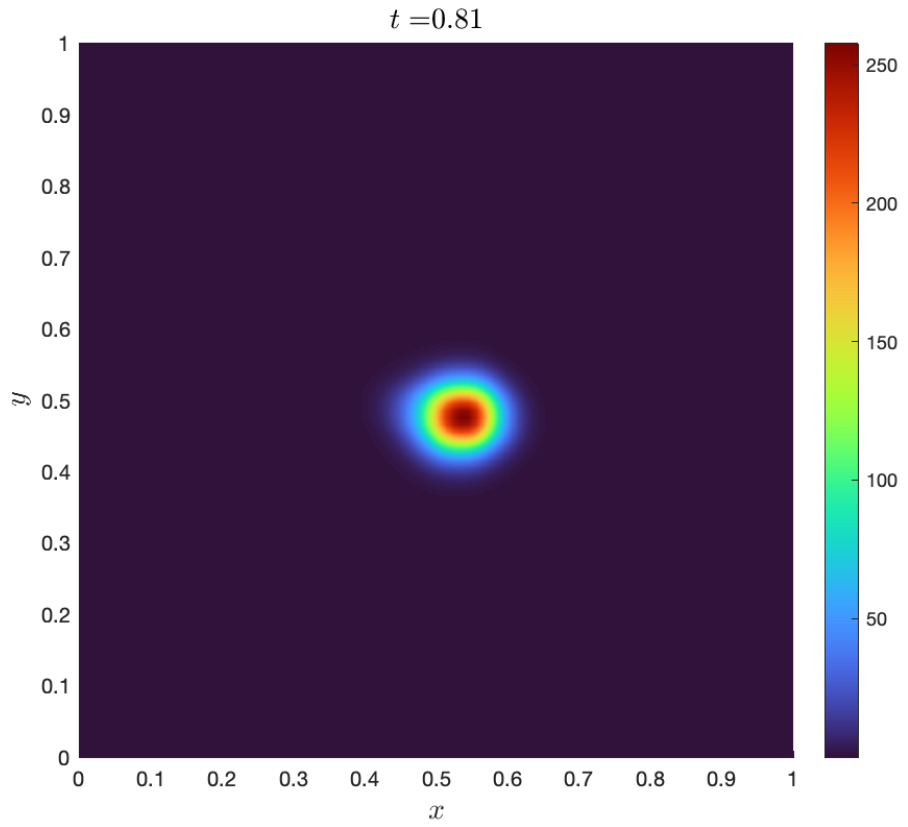}}
\caption{Numerical simulations of \eqref{eq:NP_Eulero} with $\mathcal{I}=0$, for two different values of $\eta$. Increased values of $\eta$ result in higher maximum density values, reducing the time at which the blow-up is observed. }
\label{fig:nopressnoint}
\end{figure}

Our focus is on the coupling with chemotaxis, hence we perform numerical simulations, varying the value of the parameters $\eta$. We observe that the chemotactic effect acts as an aggregating factor, directing the velocity field towards points of higher density. Fig. \ref{fig:nopressnoint} refers to the results for three different values of $\eta$. In particular, the snapshots of the density are acquired at the blow-up time instant. We can observe that the time required to detect a blow-up in the solution decreases with respect to $\eta$. This is reasonable, since higher values of the investigated parameter correspond to a strongest influence of the chemotactic gradient on the density dynamics.

\subsection{Isentropic Euler with chemotaxis}
\begin{figure}[ht]
\centering
\includegraphics[scale=0.3]{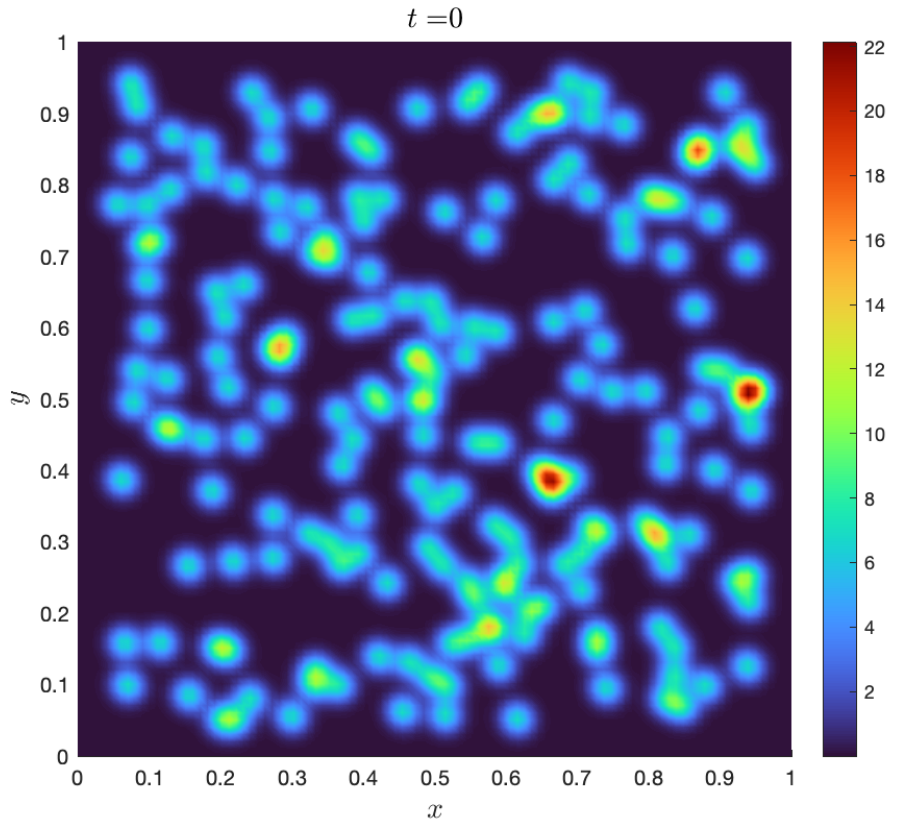}
\caption{Initial condition for $\rho$ in all simulations. }
\label{fig:density_initialcondition}
\end{figure}

Let now focus on the model in \eqref{eq:GambaPreziosi_intro}. In relevant works of the literature of vasculogenesis \cite{PhysRevLett.90.118101, ambrosi2005review, tosin2006mechanics}, the presence of the phenomenological pressure term together with the evolution of the chemical signal drives the formation of a network structure. In particular, the pressure term acts as a repulsion force, meaning that overcrowding regions are prevented. The structure of the function $P(\rho)$ to be consider is widely debated. In this simulation, we assume the pressure function introduced in
\cite{cavalli2007approximation}, i.e.
\begin{equation}
\label{pressure_qualitative_preziosi}
     P(\rho)=
    \begin{cases}
        (\rho-\rho_0)^3 & \text{if } \rho>\rho_0, \\
        0 & \text{if } \rho \le \rho_0,
    \end{cases} 
\end{equation}
with $\rho_0=4$. The presence of an activation threshold $\rho_0$ models the fact that the pressure acts as a cell repulsion at sufficiently high cell density, whereas no effect is expected at lower densities. In Fig. \ref{fig:sipress}, we observe the evolution in time of the density $\rho$, starting from the initial condition in Fig. \ref{fig:density_initialcondition}. The network formation is strictly dependent to the threshold $\rho_0$ introduced in the expression for the pressure \eqref{pressure_qualitative_preziosi}: cell migration toward other cells is strongly damped when the density reaches its threshold, since the pressure term to prevent overcrowding of cells becomes non-zero.
\subsection{Pressureless Euler with chemotaxis and non-local interactions}
We numerically investigate the non-local pressureless Euler-type model in \eqref{eq:NP_Eulero}, with the non-local interaction functions given in \eqref{eq:gamma_1}-\eqref{eq:gamma_2}. 
We consider both alignment and attraction-repulsion effects as mechanical interactions.
The structure of functions $\gamma_1$ and $\gamma_2$ are modeling choices, but depending on the ones chosen (or experimentally observed) at the microscopic scale. This is the main difference and a key factor of the model with respect to \eqref{eq:NP_Eulero}, where the choice of the pressure function is pure phenomenological.

\begin{figure}[H]
\centering
\subfigure[]{\includegraphics[scale=0.275]{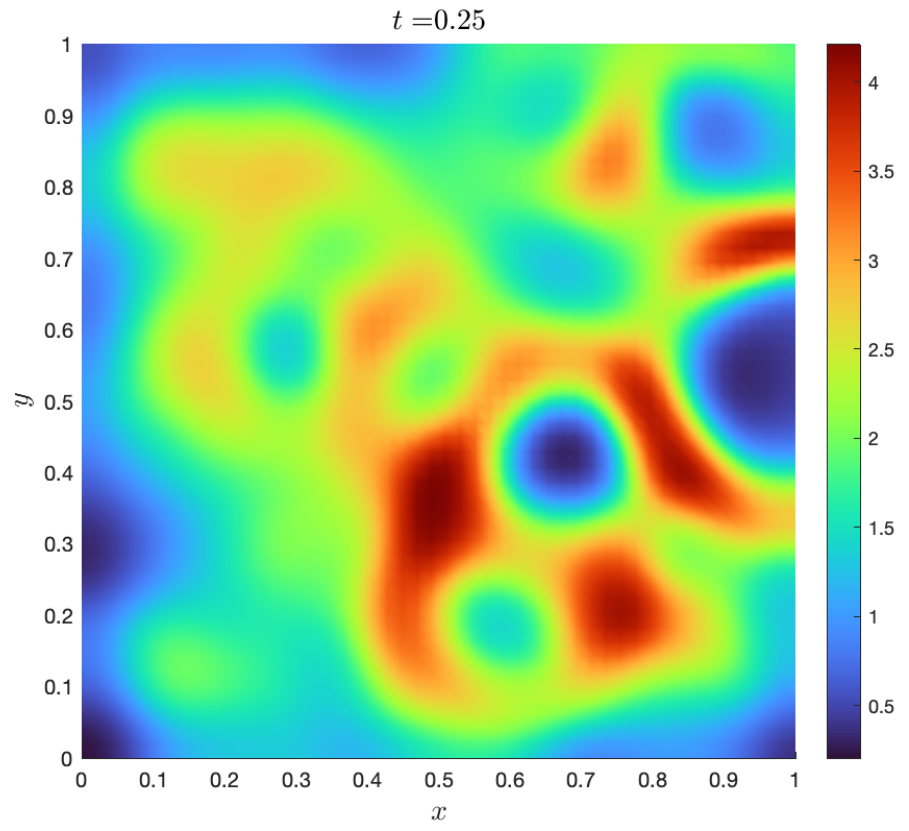}}
\subfigure[]{\includegraphics[scale=0.275]{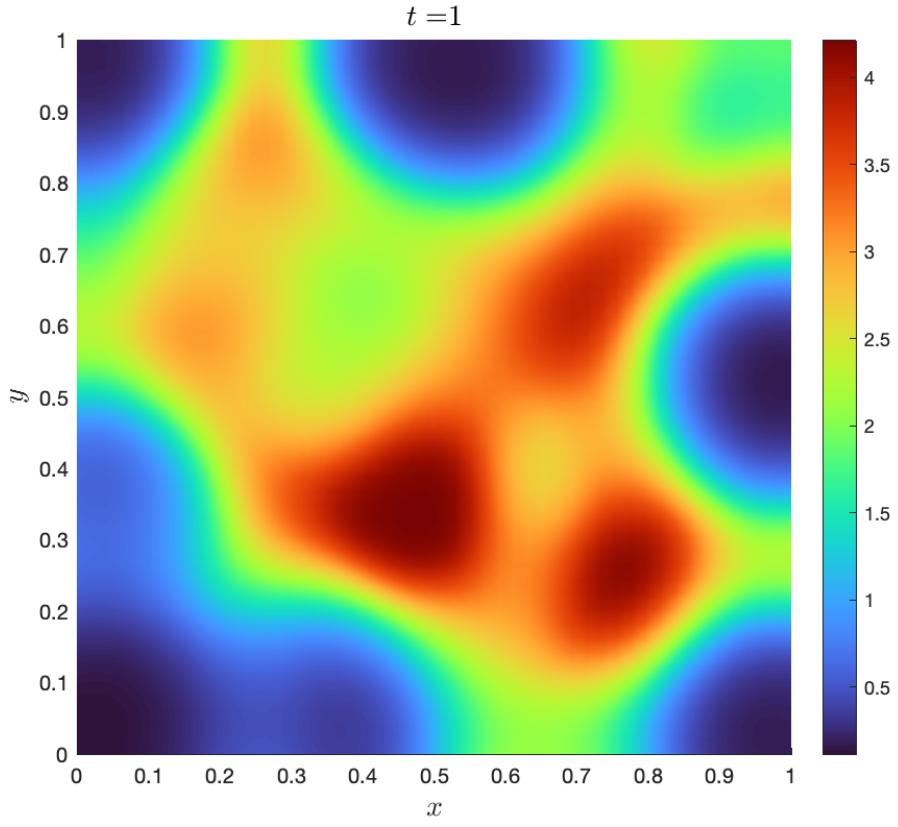}}
\subfigure[]{\includegraphics[scale=0.275]{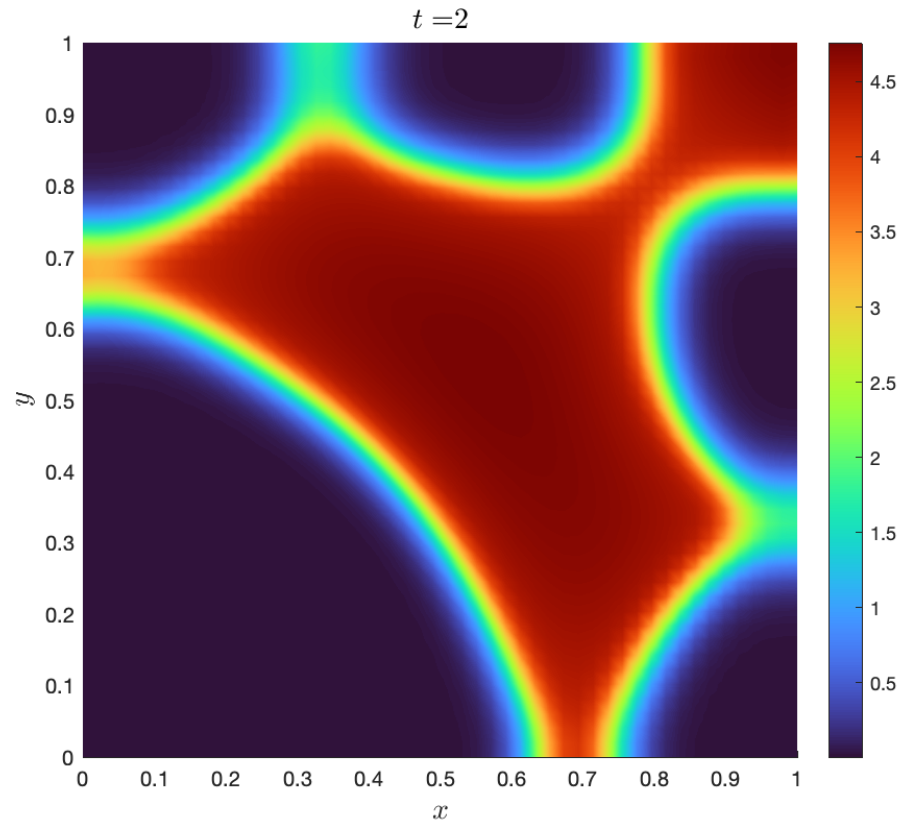}}
\caption{Numerical simulation of \eqref{eq:GambaPreziosi_intro} with $\eta=0.5$. } 
\label{fig:sipress}
\end{figure}

\begin{figure}[H]
\centering
\subfigure[]{\includegraphics[scale=0.275]{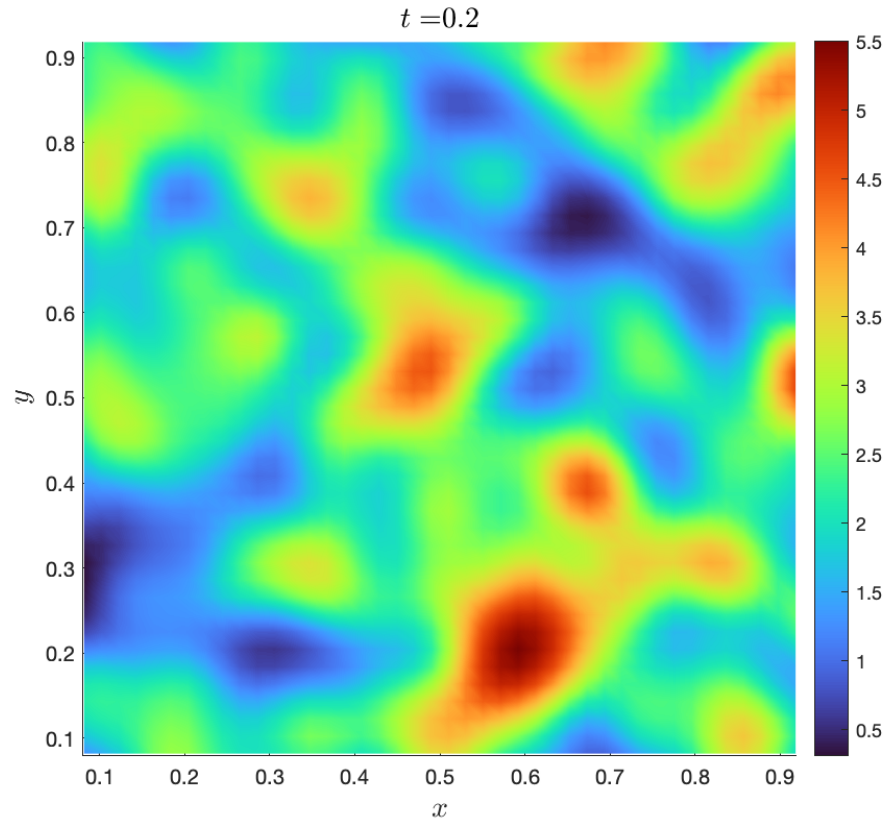}}
\subfigure[]{\includegraphics[scale=0.275]{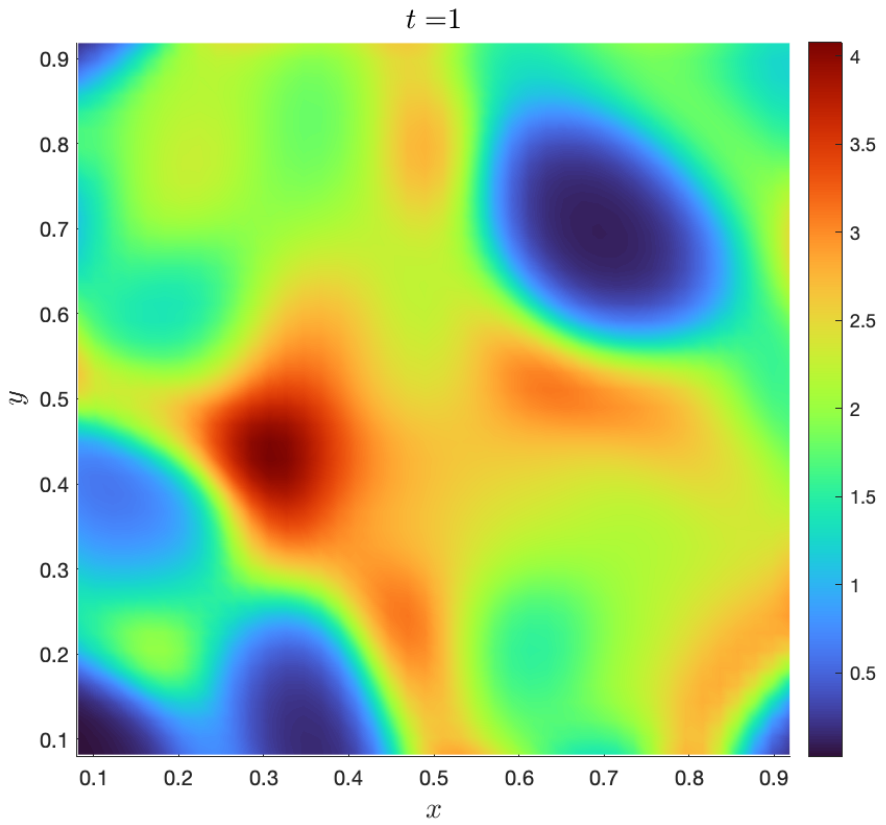}}
\subfigure[]{\includegraphics[scale=0.275]{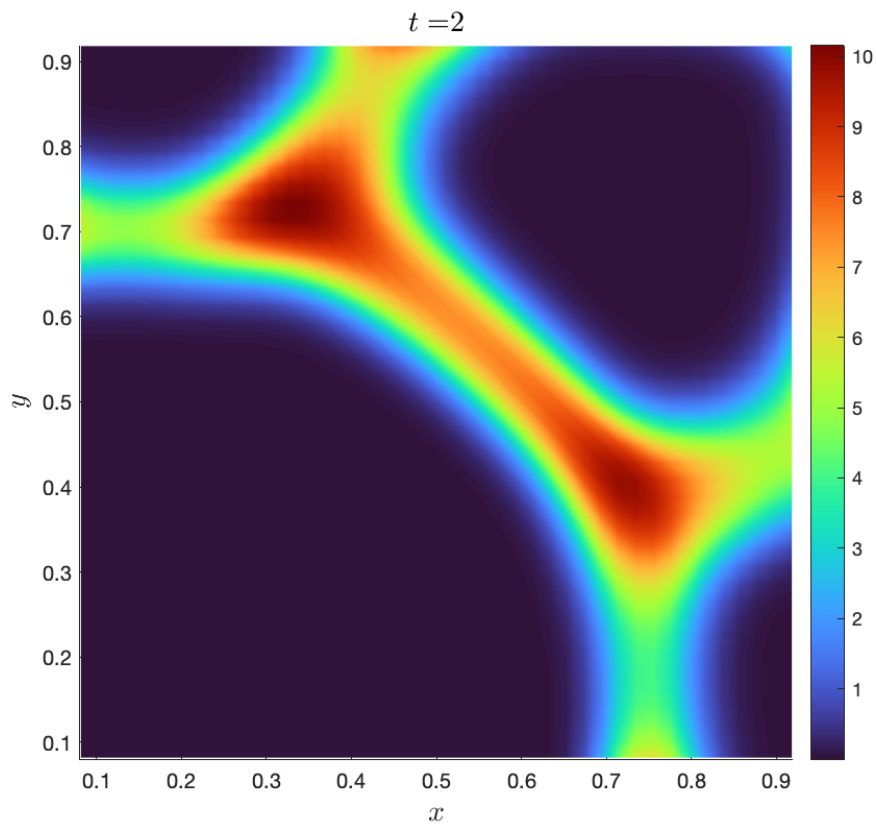}}
\caption{Numerical simulation of \eqref{eq:NP_Eulero} with $\eta=0.5$. }
\label{fig:NP_complete}
\end{figure}

Fig. \ref{fig:NP_complete} shows different snapshots of a numerical simulation of \eqref{eq:NP_Eulero}.
From a qualitative perspective, numerical evidence shows that the pressureless model is still able to reproduce network formation.
Formally, from a mathematical point of view, the pressure term is indeed replaced by the non-local interactions.
With respect to the numerical simulation in Fig. \ref{fig:sipress}, it seems that the presence of a mechanical attraction and alignment effect, in addition to chemotaxis, leads to a less dissipative behavior of the network.

\section{Parameter Estimation}\label{sec:estimation}
In this section we present the procedure followed to estimate the parameters of the model \eqref{eq:NP_Eulero_Tumor}, focusing on the key parameters that are not given in the literature. 
In particular, we validate the model on a reference scenario inspired by Cancer-on-Chip experiments, hence involving two populations of interacting cells.  Using parameter values available at the microscopic scale, obtained in \cite{bretti2021estimation, bretti2024}, we present an optimization procedure to estimate the unknown parameters, after converting the information at the microscopic scale to the macroscopic one.
In order to reduce the number of parameters to be estimated, we fix those ones known from the literature or already experimentally determined, such as cells radii, see in \cite{bretti2021estimation}.
 In our study, synthetic initial data are used to perform the calibration. We remark that, whereas the hybrid model has been already validated against real data at the microscopic scale, this is the first study concerning the calibration of the macroscopic model.

In our approach we formulate a minimization problem, aiming at finding the optimal values for the parameters of interest, minimizing the ``distance" between synthetic (or experimental) data and the numerical approximation of the density solution to \eqref{eq:NP_Eulero_Tumor}.
In fact, the density of cells represents the quantity of interest for clinicians, as well as the information which can be easily acquired through the experiments.

\subsection{Linking microscopic data to macroscopic solution}
The first step of the calibration procedure consists of generating synthetic data through the hybrid model \eqref{eq:FinalHybridModel_Tumor}, using the literature parameters available.
Synthetic data will play the role of experimental data in the formulation of the parameter estimation problem.
The kind of data we are interested in are the trajectories of cells, i.e. positions of cells at each time of the simulation.
The crucial point relies in the comparison between the microscopic information and the macroscopic counterpart.
Following the idea in \cite{braun_PhD_thesis}, we adopt a Kernel Density Estimation (KDE) approach, which is known to provide an efficient tool estimate a probability density function from a sample of observation. In our approach, we will make use of the KDE to create a macroscopic (density) counterpart starting from a discrete set of microscopic positions of cells.

\begin{defi}
\label{defi:KDE}
	Let $\mathcal{X}=\{\x_1,\dots,\x_m\}$ be a sample set of $n$-dimensional random vectors $\x_i$, distributed according to an unknown probability density function $f$. Then the multivariate kernel density estimate $\widehat{f}_{H,n}$ of $f$ is defined as
	\begin{equation}
		\widehat{f}_{H,n}(x)=\frac{1}{m} \sum_{i=1}^m K_{H,n}(x-\x_i)=\frac{1}{m \det H} \sum_{i=1}^m K(H^{-1}(x-\x_i)),
	\end{equation}
	where $x$, $\x_i \in \R^n$, the bandwidth matrix $H \in \R^{n \times n}$ and multivariate kernel function $K$.
\end{defi}

The bandwidth matrix $H$ is symmetric and positive definite, and the kernel function $K_{H,n}$ is a normalized non-negative and bounded function for all $x \in \R^n$. 

In our study $\mathcal{X}$ is the set of positions of immune cells, dealing with the two-dimensional case. For the kernel function, we choose the following multivariate Gaussian kernel function
\begin{equation}
    \displaystyle K(x)=(2\pi)^{-\frac{n}{2}} e^{-\frac{1}{2}x^Tx}.
\end{equation}

Rigorous results on the convergence of the kernel density estimate $\hat{f}_{H,n}(x)$ towards a probability density function $f$ with a bounded variation kernel function in the one-dimensional setting can be found in \cite{nadaraya1965non, chen2017tutorial}.

We point out that choice of the coefficients of the bandwidth matrix $H$ in the multi-dimensional case is fundamental to control the smoothness of the estimate. In the one-dimensional case, where the bandwidth matrix $H$ reduces to one parameter $h$, a large value of $h$ yields to an oversmoothed estimate, increasing the bias and losing some important features in the estimation. On the other hand, a small $h$ induces an undersmoothed estimate, increasing the variance and including misleading features introduced by noise.
In our 2D setting, we assume independent smoothing along the $x$ and $y$ directions, obtaining $H$ in a diagonal form. In details, we assume
\begin{equation}
    \label{bandwidth_matrix}
    H = \text{diag} ( h \sigma, h \sigma),
\end{equation} 
where $\sigma$ is the radius of the cells and $h \in \R$ is one of the parameter that will be estimated in the procedure.

\begin{figure}[h]
\centering
\subfigure[$h=1$]{\includegraphics[scale=0.25]{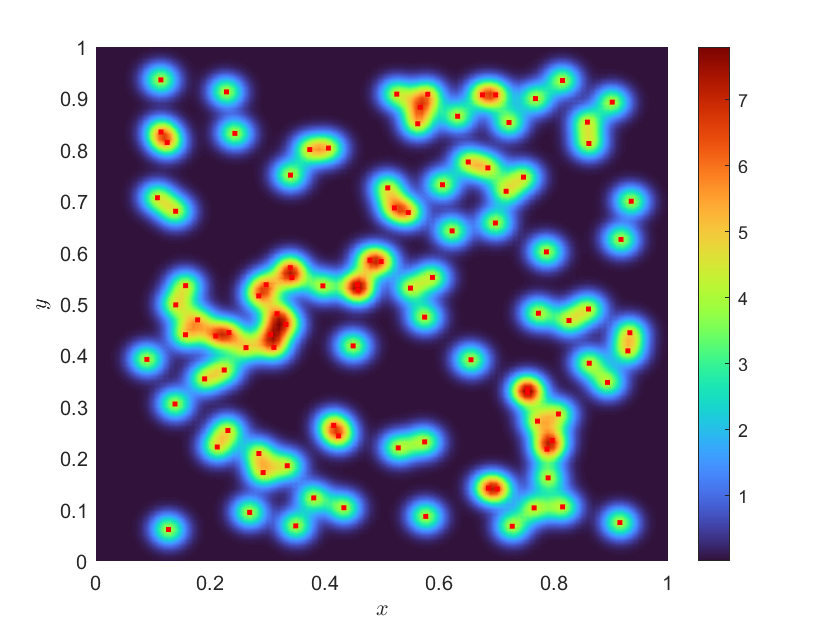}
\label{KDE_kappa=1}}
\subfigure[$h=2$]{\includegraphics[scale=0.25]{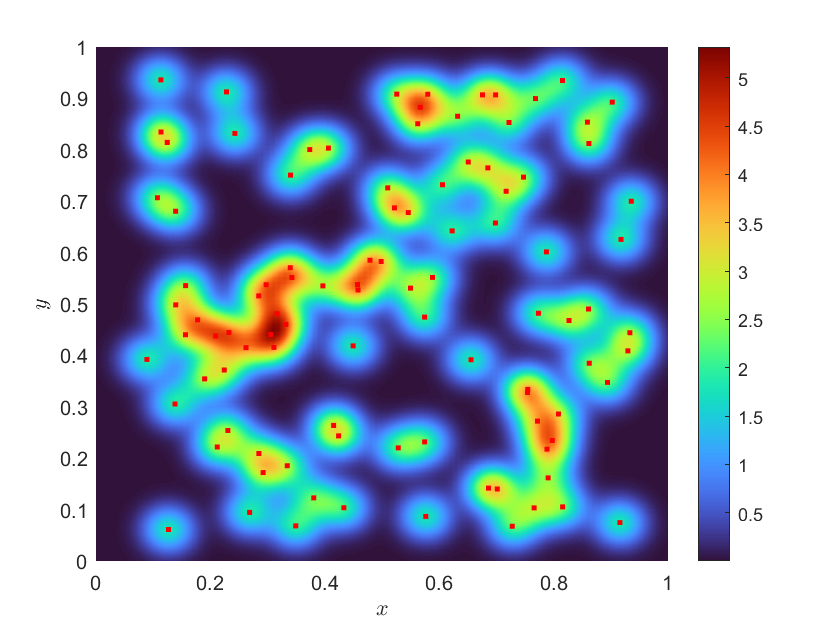}
\label{KDE_kappa=2}}
 \caption{KDE starting from a set of random microscopic positions (red dots), for two different values of $h$   in \eqref{bandwidth_matrix}. }
\label{KDE_example}
\end{figure}

In Fig. \ref{KDE_example} we show the macroscopic density obtained from a set of random microscopic positions of $m=100$ cells through the KDE approach, for two different choices of $h$.

\subsection{Minimization Problem}
To perform the parameter estimation, we define the objective function to be minimized as
\begin{equation}
    \label{objective_function_PE_KDE}
	J(\Theta)=\frac{1}{N_t} \sum_{i=1}^{N_t} \left(\frac{||\rho(t_i,\Theta)-\rho_{data}(t_i)||_2}{||\rho_{data}(t_i)||_2} \right)^2, 
\end{equation}
where $N_t$ is the number of time steps involved in the minimization process, $\rho_{data}$ is the reference solution obtained through the KDE described in the previous section and $\Theta \in \mathbb{R}^p$ is the vector of $p$ parameters to be determined.
In order to overcome possible ill-posedeness, we regularize the problem by including a Tikhonov regularization term as \cite{braun_PhD_thesis, ciavolella2024model}:
\begin{equation}
	R_\lambda(\Theta)=\lambda^2 ||\Theta-\Theta_0||_2^2,
\end{equation}
where $\Theta_0\in \mathbb{R}^p$ is a priori estimate of the parameters involved in the estimation procedure. Indeed, the presence of the regularization term allows to damp contributions from data errors and rounding errors, keeping the norm of $\lambda||(\Theta-\Theta_0)||$ of reasonable size \cite{engl1996}. In the following, we will choose small values of $\lambda$ such that the norm of $||\Theta-\Theta_0||$ will not give a significant contribution in the parameter estimation. 

The estimation can be seen as an inverse problem, aiming at solving
	\begin{equation}
		\label{eq:defi_inverse_min}
		\Theta_{opt}=\arg \min_{\Theta \in \mathcal{A}} K(\Theta),
	\end{equation}
where $\mathcal{A}$ is the set of admissible solutions and $K$ the following nonlinear regularized objective function 
	\begin{equation}
		\label{regularized_objective_function}
    K(\Theta)=J(\Theta)+R_\lambda (\Theta).
	\end{equation}

The minimization is performed using a trust region algorithm, which provides as solution the optimal value $\Theta_{opt}$. 
We assess the goodness of the parameter estimation by computing $E= J(\Theta_{opt})$, namely the relative error between synthetic data and the solution obtained with the estimated parameters.

\subsection{Numerical Results}

\begin{figure}[h]
\centering
\includegraphics[scale=0.35]{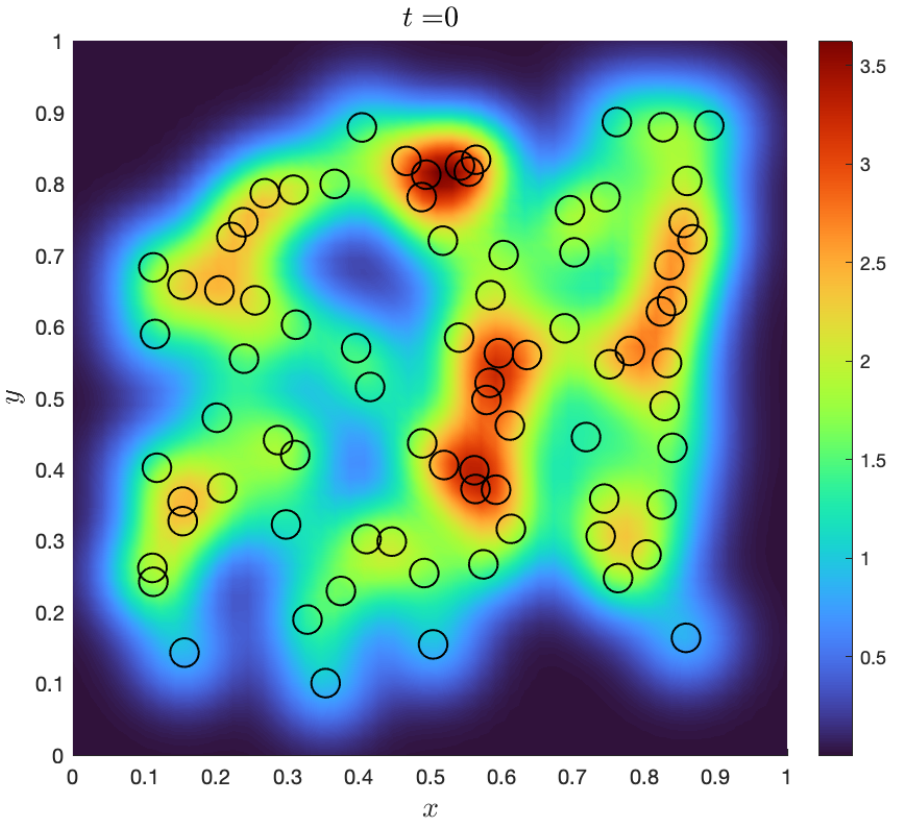}
\label{campovel_noalign}
\caption{
Plot of initial density $\rho_0$ and initial position of immune cells (black circles).
}
\label{InitialDatum_estimation}
\end{figure}

In the following we present the simulation results obtained through the parameters estimation procedure  for the macroscopic model \eqref{eq:NP_Eulero_Tumor}, accounting for two population of cells.
In order to generate synthetic data starting from the microscopic level, we run \eqref{eq:FinalHybridModel_Tumor} in the time interval $[0,1]$ considering $N=80$ immune cells, whose initial positions are chosen randomly in the domain $[0,1] \times [0,1]$ and $M=3$ tumor cells, located in $\y_1=[0.5,\, 0.7]$, $\y_2=[0.3,\, 0.3]$, $\y_3=[0.8,\, 0.5]$, respectively. Initial null velocities of cells is assumed and initial chemical concentration is chosen as in \eqref{initial_condition_tumor_chemoattractant}. 
The value of the parameters considered for the reference test are shown in Table \ref{tab::parameters}, corresponding to adimensional values of the parameters in \cite{bretti2021estimation}. 
The numerical approximation of the solution to the hybrid model, used as synthetic data in this Section, is based on an implicit Euler scheme for the ODE in \eqref{eq:FinalHybridModel_Tumor} and a Crank-Nicolson finite difference scheme for the parabolic equation, with spatial discretization performed with $\Delta x=\Delta y =0.02$. We refer to \cite{bretti2024} for a detailed description of the numerical schemes for the hybrid model. 
Starting from initial positions of immune and cancer cells, namely $\x_i$ and $\y_j$, the KDE is used to built the initial density $\rho_0$ to run the macroscopic model. We recall that cancer cells are fixed in the domain, with $\y_j(t)= \y_j$ for any $t$, hence the KDE for the cancer cell is done offline before starting the optimization procedure. In the following tests we set the bandwidth parameter value $h$ in \eqref{bandwidth_matrix} for the three cancer cells equal to $1.2$, while estimating the parameter value for the immune cells.
Fig. \ref{InitialDatum_estimation} shows the initial density $\rho_0$ plotted in the $xy$-plane, together with the positions of the immune cells, represented as circles endowed with a radius $R_{imm}$. 

The parameters to be estimated in \eqref{eq:NP_Eulero_Tumor} are  
\begin{equation}
\label{eq:Theta}
    \Theta= [\, \eta \quad \omega_{rep} \quad \omega_{adh} \quad \beta \quad \omega_{rep}^{tum} \quad h ] \in \mathbb{R}_+^6.
\end{equation}

 We recall that $\eta$ tunes the effect of the chemotactic gradient, $\omega_{rep}$, $\omega_{adh}$ are respectively the coefficients of the repulsion and the adhesion effect between immune cells \eqref{eq:gamma_2}, $\beta$ is the coefficient of the alignment term, $\omega_{rep}^{tum}$ is the coefficient of the repulsion effect between immune and cancer cells and $h$ is the bandwidth of the kernel density for the immune cells. The other parameters are fixed according to Table \ref{tab::parameters}. 

For fixed $N_t$ time steps, the numerical solution to \eqref{eq:NP_Eulero_Tumor}  is compared to the reference one when computing the value of $\Theta_{opt}$ minimizing the objective function \eqref{regularized_objective_function}. 
In the following tests, we set $N_t=5$  and $\lambda^2 = 10^{-6}$ for the regularization term.

The optimization starts from values in $\Theta_0$ given by the value of the parameters at the microscopic scale.
Then the set of admissible value $\mathcal{A}$ is chosen as $[0,\,50\Theta_0]$, taking into account physical conditions as the non-negativity of the coefficients and the fact that, by modeling construction, $R_{adh} > R_{rep}$.

\begin{table}[!h]
\footnotesize
\caption{Table of parameters values considered to generate synthetic data for Test1 - Test3.}\label{tab:foo}
\begin{center}
  \begin{tabular}{|c|c|c|} \hline
   Parameter & \bf Description & \bf Value \\ \hline
    $R_{imm}$ & Radius of ICs & $2 \cdot 10^{-2}$\\
    $R_{tum}$ & Radius of TCs & $5 \cdot 10^{-2}$\\
    $\alpha$ & Damping Coefficient & $1 \cdot 10^{\,2}$ \\ 
    $\kappa$ & Degradation Rate of $\varphi$ & $2 \cdot 10^{-1}$ \\
    $\xi$ & Production Rate of $\varphi$ & $1 \cdot 10^{\,3}$ \\
    $D$ & Diffusion Coefficient of $\varphi$ & $4.5 \cdot 10^{\,1}$ \\
    $R_{adh}$ & Radius of action of Adhesion between ICs & $6 \cdot 10^{-2}$ \\
    $R_{rep}$ & Radius of action of Repulsion between ICs & $4 \cdot 10^{-2}$ \\
    $R_{rep}^{tum}$ & Radius of action of Repulsion between ICs and TCs & $7 \cdot 10^{-2}$ \\
    $R_{align}$ & Radius of action of Alignment between ICs & $6 \cdot 10^{-2}$ \\ \hline
  \end{tabular}
\end{center}
\label{tab::parameters}
\end{table}

\paragraph{Test 1: neglecting repulsion interaction between immune and cancer cells}
In this framework, we do not consider mechanical interactions between cancer and immune cells, setting $\omega_{rep}^{tum}=0$ in  $\eqref{eq:gamma_3}$ and in $\eqref{eq:NP_Eulero_Tumor}$.
We set parameters value as in Table \ref{tab::parameters} and 
\begin{equation}
\label{Parameters_Initial_Scenario1}
    \eta=6, \quad \omega_{rep}=5 \cdot 10^2, \quad \omega_{adh}=4 , \quad \beta=2 \cdot 10^3, \quad h=1.2,
 \end{equation}
to compute the trajectories of immune cells.
Hence we run the optimization procedure, setting the components of $\Theta_0$ as the values in \eqref{Parameters_Initial_Scenario1}.
The optimization procedure gives the following optimal solution
\begin{equation}
\label{theta_opt_1}
    [71.62 \quad 6.74 \cdot 10^2  \quad 39.43 \quad 1.68 \cdot 10^3 \quad 0 \quad 2.45],
\end{equation}

with estimation error
\begin{equation}
    E= 6.73 \cdot 10^{-2}.
\end{equation}

\begin{figure}[!h]
\centering
\subfigure[ ]{\includegraphics[scale=0.3]{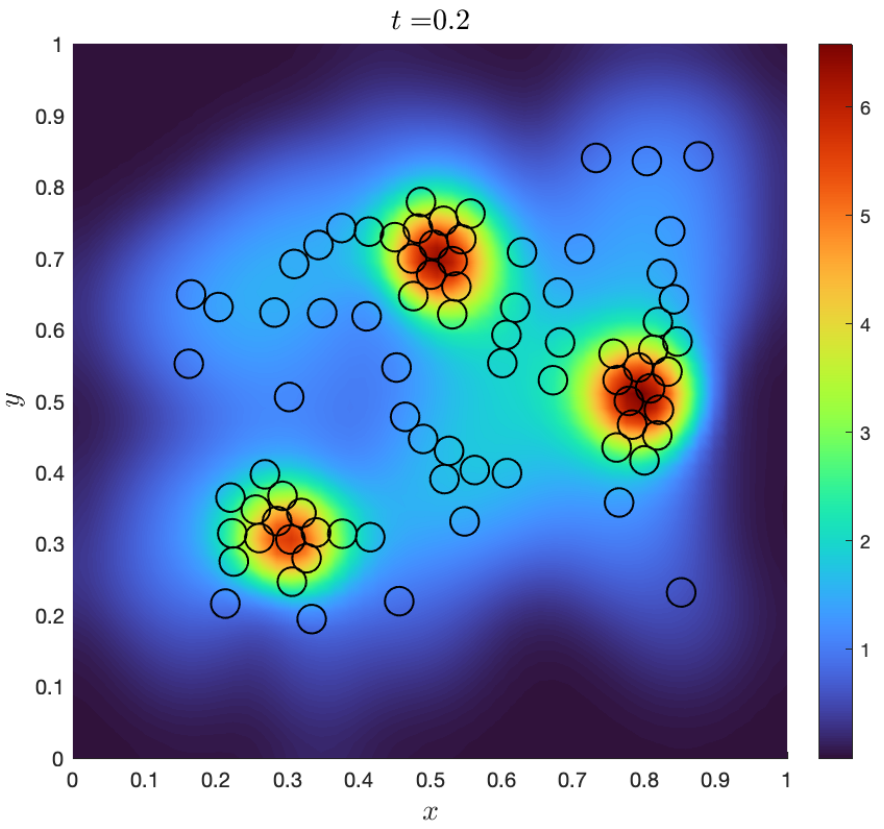}
\label{estimation_t0dot2}
}
\subfigure[ ]{\includegraphics[scale=0.3]{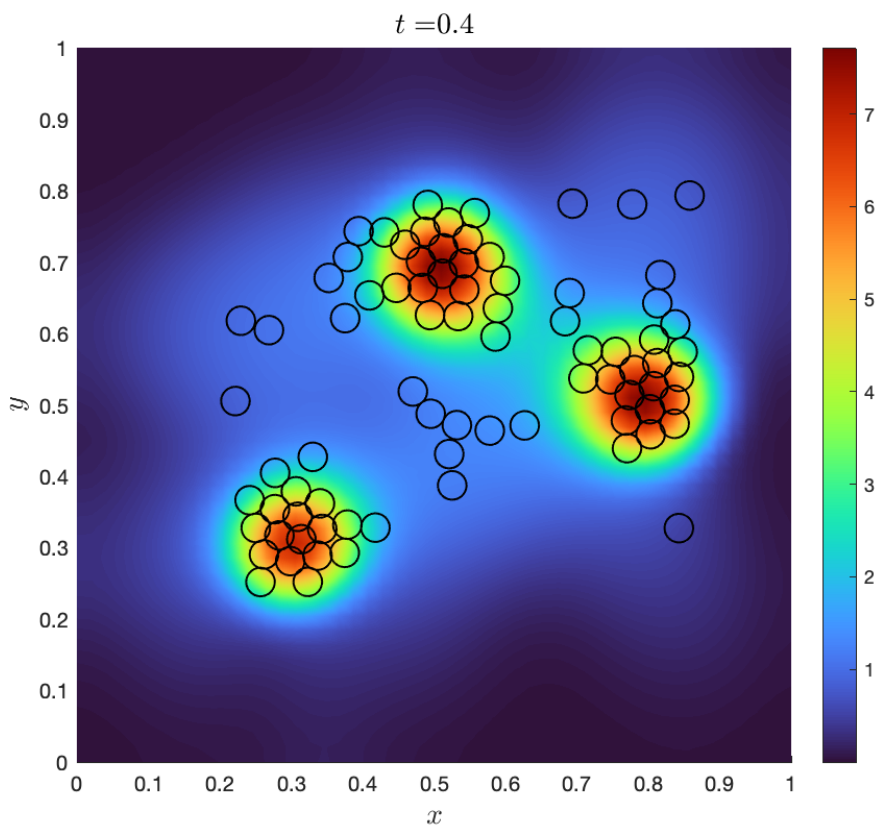}
\label{estimation_t0dot4}
}\\
\subfigure[ ]{\includegraphics[scale=0.3]{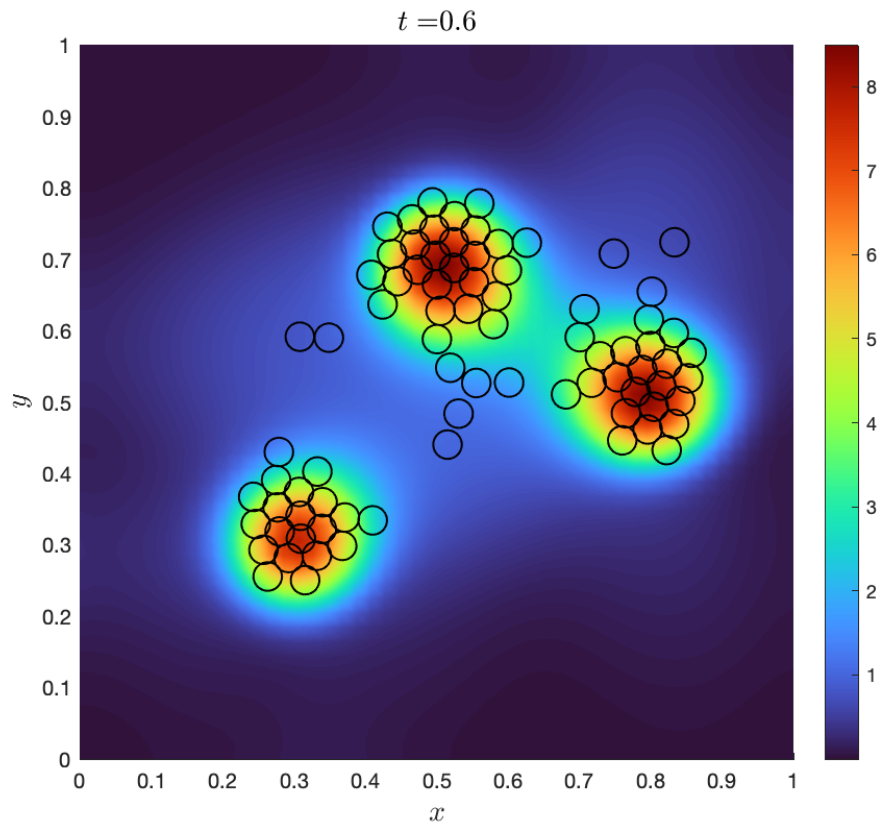}
\label{estimation_t0dot6}
}
\subfigure[ ]
{\includegraphics[scale=0.3]{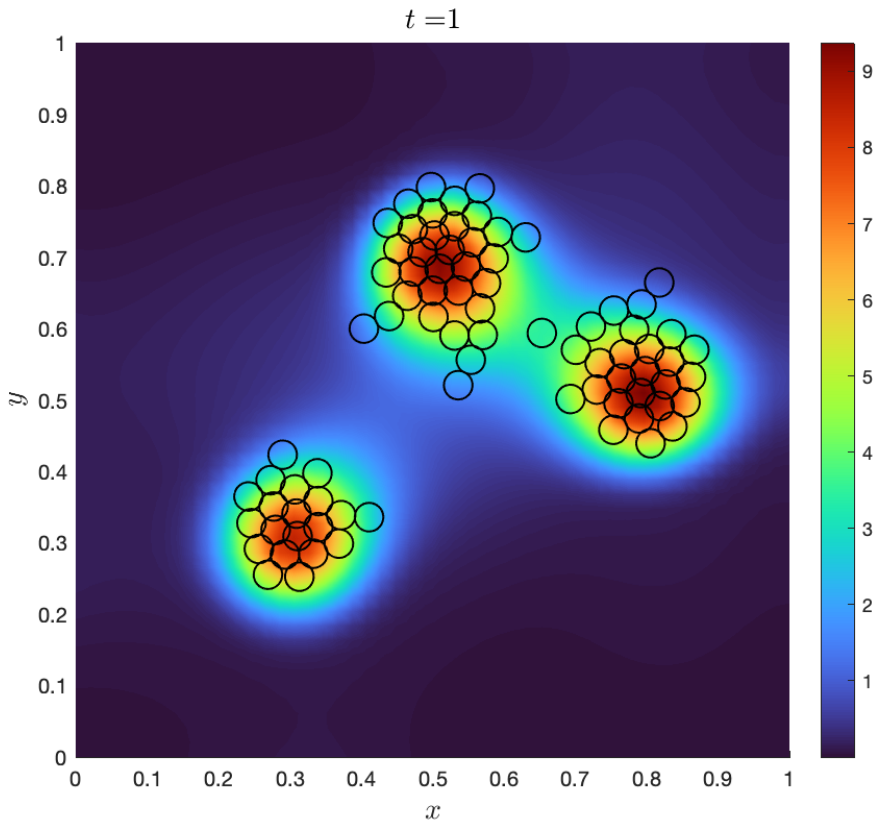}
\label{estimation_t1}
}
\caption{Numerical simulation of \eqref{eq:NP_Eulero_Tumor} with estimated parameters. (a)-(e) Plots at different times of the (macroscopic) density $\rho$ versus (microscopic) positions of the cells, represented by black circles.}
\label{Test_estimation}
\end{figure}

Fig. \ref{Test_estimation} shows five plots of the numerical simulation of $\eqref{eq:NP_Eulero_Tumor}$ with the estimated parameters at different times. Positions of the immune cells, generating the reference synthetic data, are marked as black disks. We observe that, also from a qualitative point of view, a good level of agreement is achieved between the two scales. In particular, the effect of the chemical signal produced by the three cancer cells plays a crucial role. At the microscopic scale, cells are chemically attracted by the chemical signal, and since $\omega_{rep}^{tum}=0$, they overlap cancer cells. On the microscopic scale, the density gets higher in the surrounding of $\y_1$, $\y_2$, $\y_3$.

\begin{figure}[h]
\centering

\subfigure[]{\includegraphics[scale=0.3]{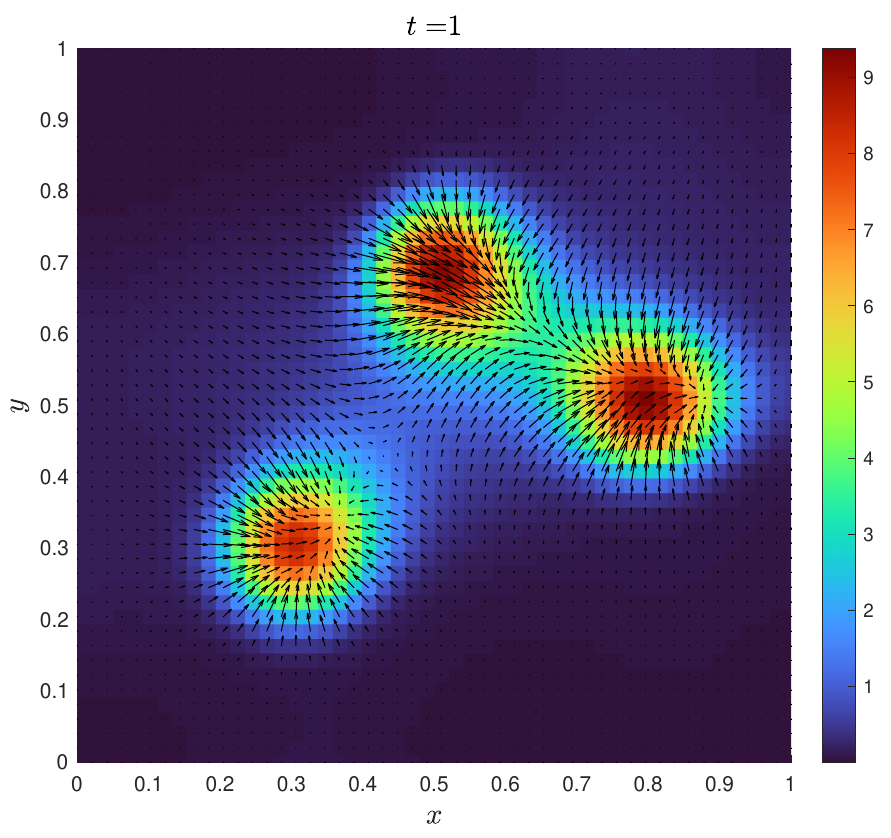}\label{campovel_sialign}
}
\subfigure[]{\includegraphics[scale=0.3]{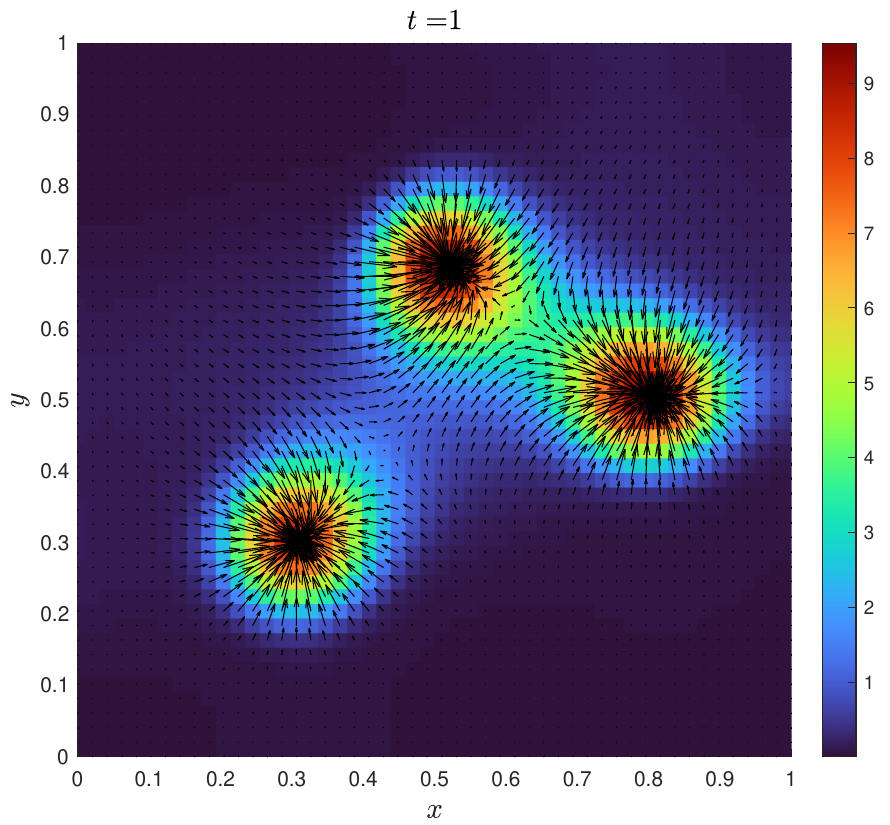}
\label{campovel_noalign}
}
\caption{Comparison at final time of numerical simulations of \eqref{eq:NP_Eulero_Tumor} with estimated parameters assuming a) $\mathcal{I}=\mathcal{I}_1 + \mathcal{I}_2 + \mathcal{I}_3$, b) $\mathcal{I}= \mathcal{I}_2 + \mathcal{I}_3$. Black arrows describing the momentum $\rho \textbf{u}$, are plotted against the plots of the density $\rho$.} 
\label{Test_campovel}
\end{figure}

\paragraph{Test 2: the effect of alignment} With respect to attraction and repulsion kind of interactions, the role of alignment in the dynamics is less evident. For this reason, we simulate the same scenario of Test 1 but neglecting alignment, hence assuming $\beta=0$ both at the microscopic level in \eqref{eq:FinalHybridModel_Tumor} and in the macroscopic one \eqref{eq:NP_Eulero_Tumor}. Consequently, we set again the components of $\Theta_0$ as the values in \eqref{Parameters_Initial_Scenario1}, with $\beta=0$.
We obtain
\begin{equation}
\label{theta_opt_2}
    \Theta_{opt} = [81.3 \quad 7.52 \cdot 10^2 \quad 39.2 \quad 0 \quad 0 \quad 2.4], 
\end{equation}
with
\begin{equation}
    E=6.01 \cdot 10^{-2}.
\end{equation}
In Fig. \ref{Test_campovel} we show the comparison of the macroscopic densities computed respectively with \eqref{theta_opt_1} and \eqref{theta_opt_2} at final time.
In the literature of hybrid models, it has already been observed that alignment plays a crucial role in the context of Cancer-on-chip modeling \cite{costanzo2020hybrid, bretti2024} at the microscopic scale. Due to the presence of alignment, the convergence of immune cells towards cancer cells less strong, allowing cells to align with others in the surroundings (Fig. a).
Neglecting alignment, immune cells are trapped by cancer cells (Fig. b), almost clustering around them.  Black arrows are used to plot the momentum $\rho \textbf{u}$ evaluated in each point of the numerical grid, highlighting the behavior of the system at the macroscopic scale.

\paragraph{Test 3: repulsion between immune and cancer cells}
In the third scenario, we consider mechanical interactions between immune and cancer cells. Despite the lack of experimental data, this could be considered a more realistic setting, in which a repulsion effect among cancer cells and immune cells occurs at a distance between the centers less than $R_{rep}^{tum}$, avoiding non feasible overlapping.
In order to generate synthetic data, we run the microscopic model \eqref{eq:FinalHybridModel_Tumor} starting with the same initial condition of the previous tests, in the same parameter setting of Table \ref{tab::parameters}, adding a non-null value for $\omega_{rep}^{tum}$.
Hence we start the optimization procedure with
\begin{equation}
\label{Parameters_Initial_Scenario3}
    \eta=6, \quad \omega_{rep}=5 \cdot 10^2, \quad \omega_{adh}=4 , \quad \beta=2 \cdot 10^3, \quad \omega_{rep}^{tum}=8.5 \cdot 10^2 \quad h=1.2,
 \end{equation}
 obtaining 

 \begin{equation}
\label{theta_opt_3}
    \Theta_{opt} = [80.2 \quad 1.32 \cdot 10^3 \quad 33.6  \quad 4.27 \cdot 10^3 \quad  2.31 \cdot 10^2 \quad 2.5], 
\end{equation}

with
\begin{equation}
    E= 9 \cdot 10^{-2}.
\end{equation}

\begin{figure}[h]
\centering
\subfigure[ ]
{\includegraphics[scale=0.3]{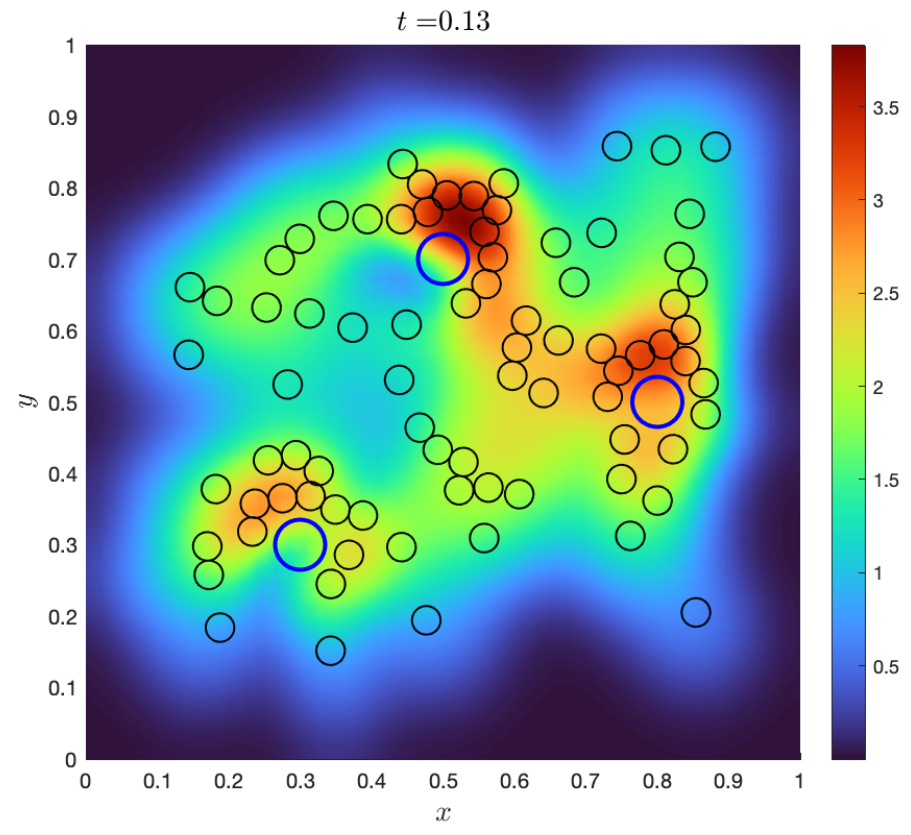}
\label{estimation_t0dot13_WREP}
}
\subfigure[ ]{\includegraphics[scale=0.3]{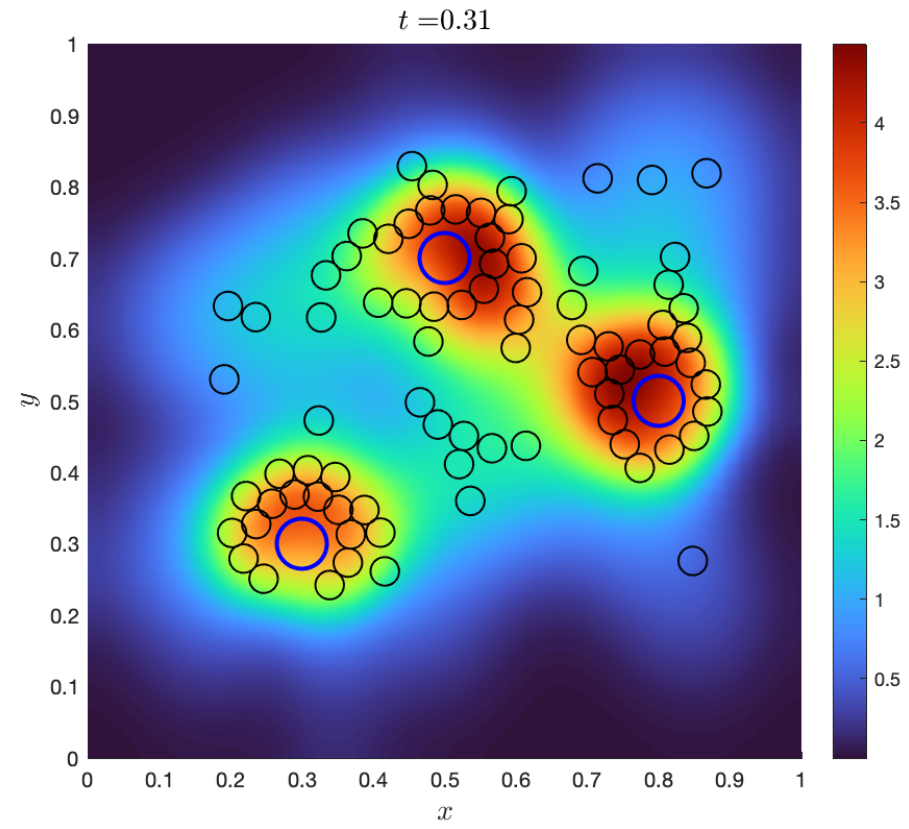}
\label{estimation_t0dot3_WREP}
}\\
\subfigure[ ]
{\includegraphics[scale=0.3]{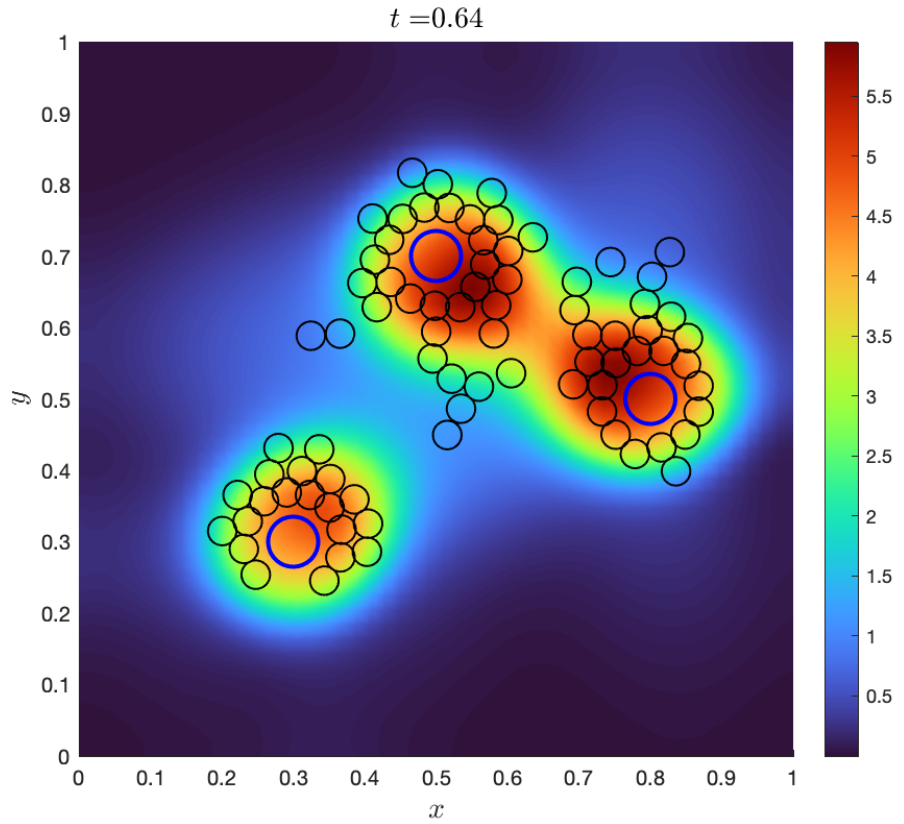}
\label{estimation_t0dot64_WREP}
}
\subfigure[ ]
{\includegraphics[scale=0.3]{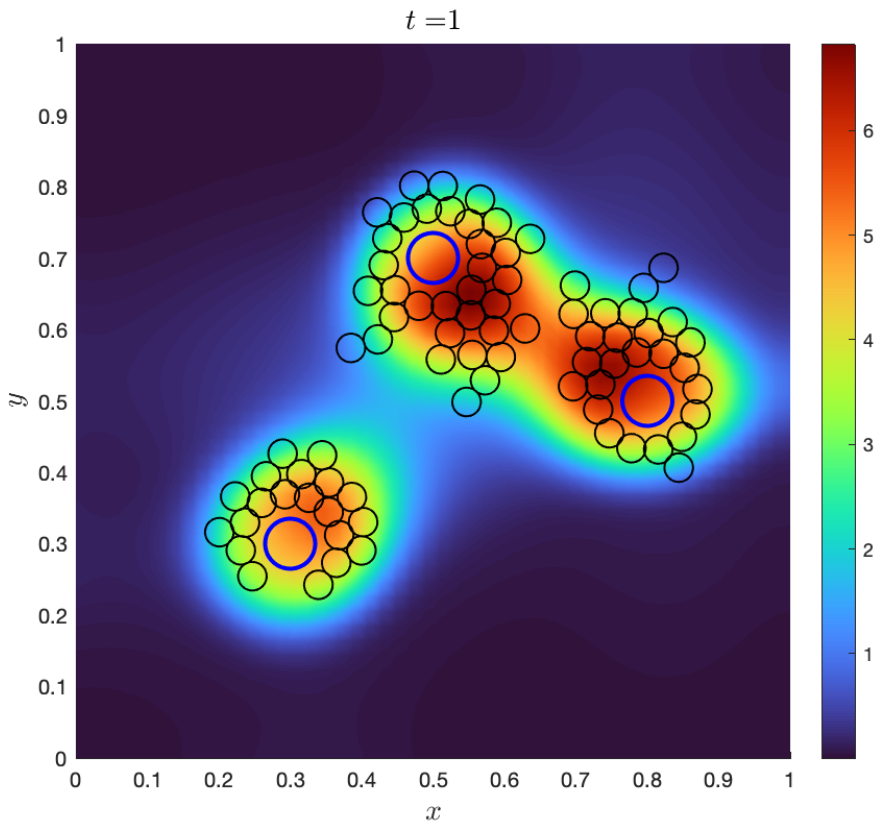}
\label{estimation_t1_WREP}
}

\caption{ Numerical simulation of  
\eqref{eq:NP_Eulero_Tumor} with estimated parameters. 
(a)-(e) Plots at different times of the (macroscopic) density $\rho$ versus (microscopic) positions of the immune and tumor cells, represented by black and blue circles, respectively. }
\label{estimation_wrepTI}
\end{figure}

Fig. \ref{estimation_wrepTI} shows different plots of the simulations with the estimated parameters. With respect to Fig. \ref{Test_estimation}, we observe the effect of the repulsion between cancer and immune cells: the density plot is less concentrated on tumor cells, plotted as blue disks endowed with radius $R_{tum}$, and immune cells do no overlap to them.

We remark that, in previous works of the literature of the related field \cite{bretti2021estimation,ciavolella2024model}, the estimation of the parameters involved in the considered models is performed against synthetic data produced by the model itself. 
Hence, it reasonable to expect a better agreement in term of relative error between the parameters used to generate synthetic data and the optimal parameters obtained through the minimization process. 
The issue we here address introduces a higher degree of difficulty, since synthetic data are generated through the hybrid model \eqref{eq:FinalHybridModel_Tumor}, which acts on a different scale compared to the model under investigation \eqref{eq:NP_Eulero_Tumor}.
The numerical results here obtained show that the non-local pressureless Euler model \eqref{eq:NP_Eulero_Tumor} represents a valid counterpart for  \eqref{eq:FinalHybridModel_Tumor} on the macroscopic scale, with the advantages of a clear link to the underlying microscopic dynamics.

\paragraph{Sensitivity Analysis}
To conclude our analysis, we evaluate the sensitivity of the model \eqref{eq:NP_Eulero_Tumor}, following the approach in \cite{ciavolella2024model, clarelli2016fluid}. Starting with the optimal values found in \eqref{theta_opt_3}, we evaluate the effects of variations of the model parameters in \eqref{eq:Theta} on the maximal density of immune cells at final time. Denoting with $\mathcal{Y}$($\Theta$)
=$\mathcal{Y}$ $(\Theta_1,...,\Theta_6)$ its value running the model with parameters given in $\Theta$, the sensitivity of the output to any parameter $\Theta_i$ is obtained as
$$
S=\frac{|\mathcal{Y}(\Theta_1, ..., \Theta_i \pm \delta,...,\Theta_6)- \mathcal{Y}(\Theta_1, ...,\Theta_6)|}{\mathcal{Y}(\Theta_1, ...,\Theta_6)} \frac{\Theta_i}{\delta} .
$$
Here we consider a $5\%$ perturbation of the parameters choosing $\delta=0.05 \Theta_{i}$.
The results are shown in Table\ref{tab::sensitivity}.
As expected, parameters $\eta$, $w_{rep}$ and $w_{rep}^{tum}$ play a significant role in the cells density. Compared to the others, parameter $w_{adh}$ tuning attraction mildly influences the output. This is reasonable, since chemotaxis already plays an adhesive role in the dynamics. The alignment parameter $\beta$ appears to have a small influence on the maximum density value, confirming the results obtained in Test2.

\begin{table}[!h]
\footnotesize
\caption{Sensitivity study with $\delta=5\%$ in the simulation set of Test3 with reference values given in \eqref{theta_opt_3}. $S_{\mathcal{Y}}$ is the sensitivity considering as model output the maximum density of immune cell at final time.}
\begin{center}
  \begin{tabular}{|c|c|} \hline
   Parameter & $S_{\mathcal{Y}}$  \\ \hline
    $\eta + \delta$ & 1.0978 \\
    $\eta - \delta$ & 1.0765 \\
    $ w_{rep} + \delta$ & 0.4238 \\
    $ w_{rep} - \delta$ & 0.4545  \\
    $ w_{adh} + \delta$ & 0.1084 \\
    $ w_{adh} - \delta$ & 0.1067\\
    $\beta + \delta$ &  0.0012 \\
    $\beta - \delta$ & 0.0012 \\
    $ w_{rep}^{tum} + \delta$ & 0.4310 \\
    $ w_{rep}^{tum} - \delta$ &  0.4843 \\
    $ h + \delta$ & 0.0124  \\
    $ h - \delta$ & 0.0106 \\ \hline
  \end{tabular}
\end{center}
\label{tab::sensitivity}
\end{table}

\newpage

\section{Conclusions}
\label{sec:conclusions}
In this paper we numerically investigate the pressureless Euler-type model rigorously derived in \cite{natalini2021mean} as macroscopic counterpart for the class of hybrid coupled models. The derivation requires strong assumptions on the initial data, which are not satisfied in real-world applications.  
Starting with a comparison with previous models of the literature, including a phenomenological pressure term, we observe the capability of this new model in reproducing typical patterns of collective cell migration.
The investigation is complemented by a quantitative calibration against data collected in Cancer-On-Chip settings. To this aim, we extend the model in different directions: including two populations of cells, with different behaviors, and including attraction-repulsion kind of interactions in addition to alignment. Extending the model to a multi-population framework enables a more realistic representation of the phenomenon, as it captures the distinct roles of various immune cell types in collective dynamics.
All the numerical simulations are performed in a multidimensional setting.
Starting from synthetic data generated at the microscopic scale, the parameter estimation procedure in three different two-dimensional scenarios shows a good level of agreement between the two scales.
The availability of new experimental data, with (macroscopic) information concerning the chemical concentration, will allow to further exploit the potential of the proposed approach. From a theoretical point of view, the non-locality of the model and the absence of a pressure term still represent challenging aspects for theoretical studies, leaving space for further investigations.

\section*{Acknowledgments}
M. Menci and T.Tenna are members of GNCS-INDAM research group, R. Natalini is member of GNAMPA-INDAM research group.
M. Menci is funded by INdAM - GNCS Project, CUP E53C24001950001, entitled ``MODA: Integrating MOdel-based and DAta-Driven Methods for Multiscale Biological Systems''.
T. Tenna is funded by the European Union’s Horizon Europe research and innovation programme under the Marie Skłodowska-Curie Doctoral Network DataHyking (Grant No.101072546).


\appendix
\section{Numerical Scheme for macroscopic models}
\label{Appendix_Numerical_Scheme}
Let us consider the Euler-type systems coupled with chemotaxis studied in this paper in a two-dimensional bounded domain $\Omega \subseteq \R^2$, rewritten in the following compact form
\begin{equation}
\label{eq:Euler_Appendix}
    \begin{cases}
    \partial_t w + \partial_{x_1} A_1(w) + \partial_{x_2} A_2(w) = F(w, \varphi),\\
    \partial_t \varphi = D \Delta \varphi + a \rho - b \varphi,
    \end{cases}
\end{equation}
where $w=(\rho, \rho u_1, \rho u_2)$ and
\begin{equation}
	A_1(w)=\begin{pmatrix}
		\rho u_1 \\
		\rho u_1^2+ \eps P(\rho)\\
		\rho u_1 u_2
	\end{pmatrix}, \qquad 
	A_2(w)=\begin{pmatrix}
		\rho u_2 \\
		\rho u_1 u_2\\
		\rho u_2^2+ \eps P(\rho)
	\end{pmatrix}.		
\end{equation}
and
\begin{equation}
    F(w, \varphi)= \begin{pmatrix}
		0 \\[7pt]
		(1-\eps) \rho\int \gamma(\cdot - y, u_1(\cdot)-u_1(y))\rho(y) dy + \eta \rho \partial_{x_1} \varphi - \alpha \rho u_1 \\[7pt]
        (1-\eps) \rho\int \gamma(\cdot - y, u_2(\cdot)-u_2(y))\rho(y) dy + \eta \rho \partial_{x_2} \varphi - \alpha \rho u_2
	\end{pmatrix}.
\end{equation}
The choice $\eps=1$ leads to \eqref{eq:GambaPreziosi_intro}, whereas the choice $\eps=0$ leads to the non-local pressureless Euler type model \eqref{eq:NP_Eulero}.\\
\paragraph{Numerical Approximation of the Euler-type system}
In order to numerically solve \eqref{eq:GambaPreziosi_intro} and \eqref{eq:NP_Eulero}, we consider a class of numerical schemes based on a discrete kinetic approximation for multidimensional systems of conservation laws, following \cite{aregba2000discrete}. This approach enables us to approximate a general quasi-linear conservation laws of the form
\begin{equation}
    \partial_t w + \nabla_x A(w)=0,
\end{equation}
by a sequence of semi-linear hyperbolic systems, through a BGK relaxation equation, leading to
\begin{equation}
	\label{semi-linear general}
	\partial_t f^\eps + \sum_{j=1}^d \Lambda_j \partial_{x_j} f^\eps=\frac{1}{\eps} \Big(M(w^\eps)-f^\eps\Big),
\end{equation}
with Cauchy initial condition $f^\eps(x,0)=f_0^\eps(x)$
and by setting
\begin{equation}
	\label{u_eps_BGK}
	u^\eps:=\displaystyle \sum_{i=1}^l f_i^\eps.
\end{equation}
Here $\eps$ is a positive number, $\Lambda_j$ are real diagonal matrices $l \times l$, $f^\eps=(f_1^\eps,\dots,f_l^\eps)$ and $\mathcal{M}: \R^k \to \R^l$ is a Lipschitz continuous function.\\
The numerical schemes presented in \cite{aregba2000discrete} are constructed by splitting \eqref{semi-linear general} into a homogeneous linear part and an ordinary differential system, which is exactly solved thanks to the particular structure of the source term.

The main advantage of this approximation, especially in the multidimensional case, is the possibility of avoiding the resolution of local Riemann problems in the design of numerical schemes. In this framework, the scalar and the system case are treated in the same way at numerical level and, since all the approximating problems are in diagonal form, the numerical scheme is easily implemented.

Here we detail the construction of the scheme for our model \eqref{eq:Euler_Appendix}.
Let us take $\Omega=[0,L] \times [0,L]$ and let use denote $\Delta x$ the space step. We consider the discretization points $\x_\alpha=(\alpha_1 \Delta x, \alpha_2 \Delta x)$, $0 \leq \alpha_i \leq N_x+1$. We denote the time step by $\Delta t$ and the approximation of a function $f$ at time $t_n=n\Delta t$ as $f^n$.
Let us introduce a 5-velocities relaxation scheme with
\begin{equation}
	\label{eq:vel_max}
	\lambda_1=\lambda(1,0), \quad \lambda_2=\lambda(0,1), \quad \lambda_3=\lambda(-1,0), \quad \lambda_4=\lambda(0,-1), \quad 
	\lambda_5=(0,0),
\end{equation}
for some $\lambda >0$. The time and space steps have to satisfy the CFL stability condition, setting $\Delta t = 0.9 \Delta x / \lambda$.\\ 
Now we introduce the corresponding Maxwellians $M_i(w) \in \R^3$, $i=1,\dots,5$
\begin{equation}
	M_i(w)=a_i w+ b_{i1} A_1(w) + b_{i2} A_2(w),
\end{equation}
where
\begin{equation}
	\begin{gathered}
		a_1=\cdots=a_4=a, a_5=1-4a, \quad a<\frac{1}{4},\\
		b_{11}=b_{22}=-b_{31}=-b_{42}= \frac{1}{2\lambda}, \quad b_{ij}=0 \quad \text{otherwise},
	\end{gathered}
\end{equation}
in order to obtain consistent conditions for the convergence of the kinetic solution to the macroscopic one.
Let us denote by $w^{n,\alpha}$ the approximation of $w$ at the point $\x_\alpha \in \R^2$ and at time $t_n$. We set the discretization  of Maxwellians as
\begin{equation}
	f_i^{n,\alpha}=M_i(w^{n,\alpha}), \quad \text{for } i=1,\dots,5.
\end{equation}
We evolve each of the functions $f_i$, $1 \leq i \leq 5$, in time through an upwind scheme by following the velocity $\lambda_i$, as
\begin{equation}
\label{eq:f_1/2}
    f^{n+\frac{1}{2},\alpha} = H^{up}(\Delta t) (f_i^{n,\alpha}),
\end{equation}
where
\begin{multline}
	H^{up}(\Delta t)(f_i^{n,\alpha})=f_i^{n,\alpha} -\mu \sum_{j=1}^2 \lambda_{ij}(f_i^{n,\alpha+e_j}-f_i^{n,\alpha-e_j})\\+\mu \sum_{j=1}^2 |\lambda_{ij}|(f_i^{n,\alpha+e_j}-2f_i^{n,\alpha}+f_i^{n,\alpha-e_j}),
\end{multline}
and $\mu=\Delta t/ (2 \Delta x)$. This upwind scheme could be improved through a flux limiter approach, in order to reduce the diffusion and obtain a second order approximation far from oscillations. Therefore, the new scheme reads as \eqref{eq:f_1/2}, where $H^{up}$ is replaced by
\begin{multline}
    H^{S} (\Delta t)= H^{up} (\Delta t) - \frac{1}{2} \mu  \sum_{j=1}^2 \lambda_{ij} \left(\text{sign}(\lambda_{ij}) - \mu \lambda_{ij} \right) \\ \left( \phi_{i,j+\frac{1}{2}}\Delta_{j+\frac{1}{2}} f_i^{n,\alpha} -  \phi_{i,j-\frac{1}{2}} \Delta_{j-\frac{1}{2}} f_i^{n,\alpha}\right),
\end{multline}
where we have introduced the following quantities 
\begin{equation}
    \phi_{i,j+\frac{1}{2}} := \phi \left(\frac{\Delta_{j-\frac{1}{2}} f_i^{n,\alpha}} {\Delta_{j+\frac{1}{2}} f_i^{n,\alpha}} \right), \qquad \Delta_{j+\frac{1}{2}} f_i^{n,\alpha}= f_i^{n,\alpha+e_j} - f_i^{n,\alpha},
\end{equation}
in which $\phi$ is the minmod function defined as $\phi (r) = \max \{0, \min\{r,1\} \}$.

Finally, we just conclude by setting
\begin{equation}
    w^{n+1}=\sum_{i=1}^5 f_i^{n+\frac{1}{2}} + F(w^{n}).
\end{equation}
The $F$ is the source term approximated in explicit as 
\begin{equation}
    F(w^n)=\begin{pmatrix}
        0\\[7pt]
        \rho^n \mathcal{I}^n + \eta \rho^n D_{x_1} \varphi^n - \alpha \rho^n u_1^n\\[7pt]
        \rho^n \mathcal{I}^n + \eta \rho^n D_{x_2} \varphi^n - \alpha \rho^n u_2^n
    \end{pmatrix},
\end{equation}
where $\mathcal{I}^n$ is the approximated value of the integral term obtained through a quadrature formula, while $D_{x_i}$ is a three points discretization of the operator $\partial_{x_i} \varphi$.\\

\paragraph{Numerical Approximation of the Parabolic Equation}
For the discretization of the parabolic equation of the chemoattractant in \eqref{eq:Euler_Appendix}, we apply a classical central finite difference scheme with 5-point stencil for the Laplacian and the $\theta$-method scheme in time. \\
The discretization of the Laplacian for a function $\varphi$ becomes:
\begin{equation}
	\Delta \varphi (x_k,y_m)=\frac{\varphi_{k+1,m}-2\varphi_{k,m}+\varphi_{k-1,m}}{\Delta x^2}+\frac{\varphi_{k,m+1}-2\varphi_{k,m}+\varphi_{k,m-1}}{\Delta y^2},
\end{equation}
while the $\theta$-method for the discretization in time leads to the final scheme for the parabolic equation:
\begin{equation*}
	\begin{aligned}
		\frac{\varphi_{k,m}^{n+1}-\varphi_{k,m}^n}{\Delta t} = &D\theta \Big(\frac{\varphi_{k+1,m}^{n+1}-2\varphi_{k,m}^{n+1}+\varphi_{k-1,m}^{n+1}}{\Delta x^2}+\frac{\varphi_{k,m+1}^{n+1}-2\varphi_{k,m}^{n+1}+\varphi_{k,m-1}^{n+1}}{\Delta y^2} \Big) \\[7pt]
		+&D (1-\theta) \Big(\frac{\varphi_{k+1,m}^n-2\varphi_{k,m}^n+\varphi_{k-1,m}^n}{\Delta x^2}+\frac{\varphi_{k,m+1}^n-2\varphi_{k,m}^n+\varphi_{k,m-1}^n}{\Delta y^2} \Big) \\[7pt]
		+&\theta s_{k,m}^{n+1}+(1-\theta) s_{k,m}^n,\\[10pt]
	\end{aligned}
\end{equation*}
where $s_{k,m}^n=s(x_k,y_m,t_n)$ is an external heat source. The numerical scheme could be complemented with suitable boundary conditions.\\
The advantage is achieved by choosing $\theta=\frac{1}{2}$, usually known as \textit{Crank-Nicolson scheme}, which increases the approximation order in time to $2$, benefiting from unconditional stability.


\bibliographystyle{plain}
\bibliography{references}

\begin{thebibliography}{10}

\bibitem{albi2018pressureless}
Giacomo Albi, Young-Pil Choi, and Axel-Stefan Haeck.
\newblock {Pressureless Euler alignment system with control}.
\newblock {\em Mathematical Models and Methods in Applied Sciences},
  28(09):1635--1664, 2018.

\bibitem{ambrosi2005review}
D~Ambrosi, F~Bussolino, and L~Preziosi.
\newblock A review of vasculogenesis models.
\newblock {\em Journal of Theoretical Medicine}, 6(1):1--19, 2005.

\bibitem{aregba2000discrete}
Denise Aregba-Driollet and Roberto Natalini.
\newblock Discrete kinetic schemes for multidimensional systems of conservation
  laws.
\newblock {\em SIAM Journal on Numerical Analysis}, 37(6):1973--2004, 2000.

\bibitem{bouchutjames1999}
Fran{\c{c}}ois Bouchut and Fran{\c{c}}ois James.
\newblock Duality solutions for pressureless gases, monotone scalar
  conservation laws, and uniqueness.
\newblock {\em Commun. Partial Differ. Equations}, 24(11-12):2173--2189, 1999.

\bibitem{bouchutjinli2003}
Fran{\c{c}}ois Bouchut, Shi Jin, and Xiantao Li.
\newblock Numerical approximations of pressureless and isothermal gas dynamics.
\newblock {\em SIAM J. Numer. Anal.}, 41(1):135--158, 2003.

\bibitem{braun_PhD_thesis}
Elishan~Christian Braun.
\newblock {\em {Organs-On-Chips: mathematical modelling and parameter
  estimation}}.
\newblock {PhD Thesis}, {Universit{\`a} degli studi Roma Tre}, 2021.

\bibitem{breniergrenier1998}
Yann Brenier and Emmanuel Grenier.
\newblock Sticky particles and scalar conservation laws.
\newblock {\em SIAM J. Numer. Anal.}, 35(6):2317--2328, 1998.

\bibitem{bretti2024}
Gabriella Bretti, Elio Campanile, Marta Menci, and Roberto Natalini.
\newblock A scenario-based study on hybrid {PDE-ODE} model for cancer-on-chip
  experiment.
\newblock In {\em Problems in Mathematical Biophysics: A Volume in Memory of
  Alberto Gandolfi}, pages 37--57. Springer, 2024.

\bibitem{bretti2021estimation}
Gabriella Bretti, Adele De~Ninno, Roberto Natalini, Daniele Peri, and Nicole
  Roselli.
\newblock {Estimation Algorithm for a Hybrid PDE–ODE Model Inspired by
  Immunocompetent Cancer-on-Chip Experiment}.
\newblock {\em Axioms}, 10(4):243, 2021.

\bibitem{carrillo2017review}
Jos{\'e}~A Carrillo, Young-Pil Choi, and Sergio~P Perez.
\newblock A review on attractive--repulsive hydrodynamics for consensus in
  collective behavior.
\newblock {\em Active Particles, Volume 1: Advances in Theory, Models, and
  Applications}, pages 259--298, 2017.

\bibitem{carrillo2016critical}
Jos{\'e}~A Carrillo, Young-Pil Choi, Eitan Tadmor, and Changhui Tan.
\newblock {Critical thresholds in 1D Euler equations with non-local forces}.
\newblock {\em Mathematical Models and Methods in Applied Sciences},
  26(01):185--206, 2016.

\bibitem{carrillo2013new}
Jos{\'e}~A Carrillo, Stephan Martin, and Vladislav Panferov.
\newblock A new interaction potential for swarming models.
\newblock {\em Physica D: Nonlinear Phenomena}, 260:112--126, 2013.

\bibitem{cavalli2007approximation}
Fausto Cavalli, Andrea Gamba, Giovanni Naldi, and Matteo Semplice.
\newblock {Approximation of 2D and 3D models of chemotactic cell movement in
  vasculogenesis}.
\newblock {\em Math Everywhere: Deterministic and Stochastic Modelling in
  Biomedicine, Economics and Industry. Dedicated to the 60th Birthday of
  Vincenzo Capasso}, pages 179--191, 2007.

\bibitem{chen2017tutorial}
Yen-Chi Chen.
\newblock A tutorial on kernel density estimation and recent advances.
\newblock {\em Biostatistics \& Epidemiology}, 1(1):161--187, 2017.

\bibitem{ciavolella2024model}
Giorgia Ciavolella, Nathalie Ferrand, Mich{\`e}le Sabbah, Beno{\^\i}t Perthame,
  and Roberto Natalini.
\newblock A model for membrane degradation using a gelatin invadopodia assay.
\newblock {\em Bulletin of Mathematical Biology}, 86(3):30, 2024.

\bibitem{clarelli2016fluid}
Fabrizio Clarelli, Cristiana Di~Russo, Roberto Natalini, and Magali Ribot.
\newblock A fluid dynamics multidimensional model of biofilm growth: stability,
  influence of environment and sensitivity.
\newblock {\em Mathematical medicine and biology: a journal of the IMA},
  33(4):371--395, 2016.

\bibitem{couzin2005effective}
Iain~D Couzin, Jens Krause, Nigel~R Franks, and Simon~A Levin.
\newblock Effective leadership and decision-making in animal groups on the
  move.
\newblock {\em Nature}, 433(7025):513--516, 2005.

\bibitem{cucker2007emergent}
Felipe Cucker and Steve Smale.
\newblock Emergent behavior in flocks.
\newblock {\em IEEE Transactions on Automatic Control}, 52(5):852--862, 2007.

\bibitem{cucker2007mathematics}
Felipe Cucker and Steve Smale.
\newblock On the mathematics of emergence.
\newblock {\em Japanese Journal of Mathematics}, 2:197--227, 2007.

\bibitem{di2018discrete}
Ezio Di~Costanzo, Alessandro Giacomello, Elisa Messina, Roberto Natalini,
  Giuseppe Pontrelli, Fabrizio Rossi, Robert Smits, and Monika Twarogowska.
\newblock {A discrete in continuous mathematical model of cardiac progenitor
  cells formation and growth as spheroid clusters (Cardiospheres)}.
\newblock {\em Mathematical Medicine and Biology: a Journal of the IMA},
  35(1):121--144, 2018.

\bibitem{costanzo2020hybrid}
Ezio Di~Costanzo, Marta Menci, Eleonora Messina, Roberto Natalini, and Antonia
  Vecchio.
\newblock A hybrid model of collective motion of discrete particles under
  alignment and continuum chemotaxis.
\newblock {\em Discrete \& Continuous Dynamical Systems-Series B},
  25(1):443--472, 2020.

\bibitem{di2015hybrid}
Ezio Di~Costanzo, Roberto Natalini, and Luigi Preziosi.
\newblock A hybrid mathematical model for self-organizing cell migration in the
  zebrafish lateral line.
\newblock {\em Journal of Mathematical Biology}, 71(1):171--214, 2015.

\bibitem{d2006self}
Maria~R D’Orsogna, Yao-Li Chuang, Andrea~L Bertozzi, and Lincoln~S Chayes.
\newblock Self-propelled particles with soft-core interactions: patterns,
  stability, and collapse.
\newblock {\em Physical Review Letters}, 96(10):104302, 2006.

\bibitem{engl1996}
Heinz~W. Engl, Martin Hanke, and Andreas Neubauer.
\newblock {\em Regularization of inverse problems}, volume 375 of {\em Math.
  Appl., Dordr.}
\newblock Dordrecht: Kluwer Academic Publishers, 1996.

\bibitem{figalli2018rigorous}
Alessio Figalli and Moon-Jin Kang.
\newblock {A rigorous derivation from the kinetic Cucker--Smale model to the
  pressureless Euler system with nonlocal alignment}.
\newblock {\em Analysis \& PDE}, 12(3):843--866, 2018.

\bibitem{PhysRevLett.90.118101}
A.~Gamba, D.~Ambrosi, A.~Coniglio, A.~de~Candia, S.~Di~Talia, E.~Giraudo,
  G.~Serini, L.~Preziosi, and F.~Bussolino.
\newblock {Percolation, Morphogenesis, and Burgers Dynamics in Blood Vessels
  Formation}.
\newblock {\em Phys. Rev. Lett.}, 90:118101, Mar 2003.

\bibitem{hillen2009user}
Thomas Hillen and Kevin~J Painter.
\newblock {A user’s guide to PDE models for chemotaxis}.
\newblock {\em Journal of Mathematical Biology}, 58(1-2):183, 2009.

\bibitem{keller1970initiation}
Evelyn~F Keller and Lee~A Segel.
\newblock Initiation of slime mold aggregation viewed as an instability.
\newblock {\em Journal of Theoretical Biology}, 26(3):399--415, 1970.

\bibitem{keller1971model}
Evelyn~F Keller and Lee~A Segel.
\newblock Model for chemotaxis.
\newblock {\em Journal of theoretical biology}, 30(2):225--234, 1971.

\bibitem{kowalczyk2004stability}
R~Kowalczyk, A~Gamba, and L~Preziosi.
\newblock On the stability of homogeneous solutions to some aggregation models.
\newblock {\em Discrete and Continuous Dynamical Systems Series B},
  4(1):203--220, 2004.

\bibitem{menci2023microscopic}
Marta Menci, Roberto Natalini, and Thierry Paul.
\newblock Microscopic, kinetic and hydrodynamic hybrid models of collective
  motions with chemotaxis: a numerical study.
\newblock {\em Mathematics and Mechanics of Complex Systems}, 12(1):47--83,
  2023.

\bibitem{menci2019global}
Marta Menci and Marco Papi.
\newblock Global solutions for a path-dependent hybrid system of differential
  equations under parabolic signal.
\newblock {\em Nonlinear Analysis}, 184:172--192, 2019.

\bibitem{menci2022existence}
Marta Menci and Marco Papi.
\newblock Existence of solutions for hybrid systems of differential equations
  under exogenous information with discontinuous source term.
\newblock {\em Nonlinear Analysis}, 221:112885, 2022.

\bibitem{menci2023coupled}
Marta Menci, Marco Papi, Maria~Michaela Porzio, and Flavia Smarrazzo.
\newblock On a coupled hybrid system of nonlinear differential equations with a
  nonlocal concentration.
\newblock {\em Journal of Differential Equations}, 361:288--338, 2023.

\bibitem{menci2023existence}
Marta Menci, Marco Papi, and Flavia Smarrazzo.
\newblock Existence of solutions for coupled hybrid systems of differential
  equations for microscopic dynamics and local concentrations.
\newblock {\em Communications on Pure and Applied Analysis}, 22(7):2146--2168,
  2023.

\bibitem{motsch2011new}
Sebastien Motsch and Eitan Tadmor.
\newblock A new model for self-organized dynamics and its flocking behavior.
\newblock {\em Journal of Statistical Physics}, 144:923--947, 2011.

\bibitem{murray2007mathematical}
James~D Murray.
\newblock {\em {Mathematical biology: I. An introduction}}, volume~17.
\newblock Springer Science \& Business Media, 2007.

\bibitem{nadaraya1965non}
EA~Nadaraya.
\newblock On non-parametric estimates of density functions and regression
  curves.
\newblock {\em Theory of Probability \& Its Applications}, 10(1):186--190,
  1965.

\bibitem{nanjundiah1973chemotaxis}
Vidyanand Nanjundiah.
\newblock Chemotaxis, signal relaying and aggregation morphology.
\newblock {\em Journal of Theoretical Biology}, 42(1):63--105, 1973.

\bibitem{natalini1998discrete}
Roberto Natalini.
\newblock A discrete kinetic approximation of entropy solutions to
  multidimensional scalar conservation laws.
\newblock {\em Journal of Differential Equations}, 148(2):292--317, 1998.

\bibitem{natalini2020mean}
Roberto Natalini and Thierry Paul.
\newblock {On the mean field limit for Cucker-Smale models}.
\newblock {\em Discrete and Continuous Dynamical Systems-B}, 27(5):2873--2889,
  2022.

\bibitem{natalini2021mean}
Roberto Natalini and Thierry Paul.
\newblock The mean-field limit for hybrid models of collective motions with
  chemotaxis.
\newblock {\em SIAM Journal on Mathematical Analysis}, 55(2):900--928, 2023.

\bibitem{natalini2015numerical}
Roberto Natalini, Magali Ribot, and Monika Twarogowska.
\newblock A numerical comparison between degenerate parabolic and quasilinear
  hyperbolic models of cell movements under chemotaxis.
\newblock {\em Journal of Scientific Computing}, 63:654--677, 2015.

\bibitem{painter2019mathematical}
Kevin~J Painter.
\newblock Mathematical models for chemotaxis and their applications in
  self-organisation phenomena.
\newblock {\em Journal of Theoretical Biology}, 481:162--182, 2019.

\bibitem{perthame2006transport}
Beno{\^\i}t Perthame.
\newblock {\em Transport equations in biology}.
\newblock Springer Science \& Business Media, 2006.

\bibitem{scianna2013}
M.~Scianna, C.~G. Bell, and Luigi Preziosi.
\newblock A review of mathematical models for the formation of vascular
  networks.
\newblock {\em J. Theor. Biol.}, 333:174--209, 2013.

\bibitem{serini2003modeling}
Guido Serini, Davide Ambrosi, Enrico Giraudo, Andrea Gamba, Luigi Preziosi, and
  Federico Bussolino.
\newblock Modeling the early stages of vascular network assembly.
\newblock {\em The EMBO journal}, 22(8):1771--1779, 2003.

\bibitem{szabo2006phase}
Balint Szabo, GJ~Sz{\"o}ll{\"o}si, B~G{\"o}nci, Zs~Jur{\'a}nyi, David Selmeczi,
  and Tam{\'a}s Vicsek.
\newblock Phase transition in the collective migration of tissue cells:
  experiment and model.
\newblock {\em Physical Review E}, 74(6):061908, 2006.

\bibitem{tosin2006mechanics}
Andrea Tosin, Davide Ambrosi, and Luigi Preziosi.
\newblock Mechanics and chemotaxis in the morphogenesis of vascular networks.
\newblock {\em Bulletin of Mathematical Biology}, 68:1819--1836, 2006.

\end{thebibliography}

\end{document}